\documentclass{amsart}
\usepackage{amsfonts,amsmath,amssymb}
\newtheorem{theorem}{Theorem}
\newtheorem{prop}[theorem]{Proposition}
\newtheorem{cor}[theorem]{Corollary}
\theoremstyle{definition}
\newtheorem*{defi}{Definition}
\newtheorem{lemma}[theorem]{Lemma}
\newtheorem{example}[theorem]{Example}
\theoremstyle{remark}
\newtheorem{remark}[theorem]{Remark}
\numberwithin{theorem}{section} \numberwithin{equation}{section}
\def\d{\partial}
\def\t{\theta}
\def\o{\omega}
\begin{document}
\author{Oleg Chalykh$^1$, Pavel Etingof, Alexei Oblomkov}
\thanks {$^1$ On leave of absence from:  Advanced Education and Science Centre, Moscow State
University, Moscow 119899, Russia}
\address{O.C.: Department of Mathematics,
Cornell University, Ithaca, NY 14853, USA}
\address{P.E.: Department of Mathematics, MIT, 77 Mass. Ave,
Cambridge, MA 02139, USA}
\address{A.O.: Department of Mathematics, MIT, 77 Mass. Ave,
Cambridge, MA 02139, USA} \email{oleg@math.cornell.edu,
etingof@math.mit.edu, oblomkov@math.mit.edu}

\title{Generalized Lam\'e operators}
\begin{abstract}
We introduce a class of multidimensional Schr\"odinger operators
with elliptic potential which generalize the classical Lam\'e
operator to higher dimensions. One natural example is the
Calogero--Moser operator, others are related to the root systems
and their deformations. We conjecture that these operators are
algebraically integrable, which is a proper generalization of the
finite-gap property of the Lam\'e operator. Using earlier results
of Braverman, Etingof and Gaitsgory, we prove this under
additional assumption of the usual, Liouville integrability. In
particular, this proves the Chalykh--Veselov conjecture for the
elliptic Calogero--Moser problem for all root systems. We also
establish algebraic integrability in all known two-dimensional
cases. A general procedure for calculating the Bloch
eigenfunctions is explained. It is worked out in detail for two
specific examples: one is related to $B_2$ case, another one is a
certain deformation of the $A_2$ case. In these two cases we also
obtain similar results for the discrete versions of these
problems, related to the difference operators of
Macdonald--Ruijsenaars type.
\end{abstract}
\maketitle

\section{Introduction}
In this paper we consider higher-dimensional analogues of the
classical Lam\'e operator
\begin{equation}\label{lame}
L=-\frac{d^2}{dz^2}+m(m+1)\wp(z)\,,\qquad m\in\mathbb Z_+\,.
\end{equation}
Here $\wp(z)=\wp(z|1,\tau)$ is the Weierstrass $\wp$-function with
periods $1,\tau$. More generally, we are interested in
multivariable analogues of the so-called {\bf elliptic
algebro-geometric operators} $L=-d^2/dz^2+u(z)$, which appeared in
the finite-gap theory initiated in 70's by Novikov \cite{N}. This
theory provides a beautiful interplay between the spectral theory
and algebraic geometry, and the Lam\'e operator is the simplest
and best known member of this family of operators (see \cite{KN}
for a survey). Since then there have been several attempts to
generalize some parts of that theory to higher dimensions, most
notably \cite{Na1,Na2}, see also \cite{Pa,Os,Ro}. We should
stress, however, that in general this leads to differential
operators with {\it matrix coefficients}. On the other hand, in
\cite{CV90} it was suggested to consider the quantum elliptic
Calogero--Moser problem and its versions related to the root
systems \cite{OP} as natural multidimensional analogues of the
Lam\'e operator. More specifically, a conjecture from \cite{CV90}
says that for integer values of the coupling parameters the
corresponding Schr\"odinger operators are {\it algebraically
integrable} (this is a proper generalization of the properties of
the algebro-geometric operators to higher dimensions, see
\cite{VSC,BEG} and Section \ref{ai} below). For the rational and
trigonometric versions of the Calogero--Moser problem this was
proved in \cite{VSC}. Elliptic version, however, turned out to be
more difficult: until now it was known for $A_n$ case only, due to
\cite{BEG}.

One of the results of the present paper is a proof of that
conjecture of \cite{CV90} for all root systems. In fact, our
approach applies to a wider class of Schr\"odinger operators,
which is an elliptic version of the class introduced in
\cite{C98}, see also \cite{CFV}. Their singularities are the
second order poles along a set of hyperplanes satisfying some
special conditions, which encode the {\it triviality of the local
monodromy} around each of the poles. We call the corresponding
Schr\"odinger operators the {\bf generalized Lam\'e operators}.
The Calogero--Moser operators with integer coupling parameters
give particular examples of such operators. Our main result says
that for a generalized Lam\'e operator its algebraic integrability
follows from the usual, Liouville integrability (in a slightly
stronger sense). The proof uses a criterion from \cite{BEG}, based
on differential Galois theory. In dimension one we recover in this
way the main result of \cite{GW1}. For the elliptic
Calogero--Moser problem the complete integrability was proved (for
all root systems) by Cherednik \cite{Ch}, and this allows us to
prove the conjecture of \cite{CV90}.

A complete description of all generalized Lam\'e operators is an
open problem. All known (irreducible) examples in dimension $>1$
are related to the root systems and their deformations which
appeared in \cite{CFV,CV1}; we list them all in Section 4. We
conjecture that they are all algebraically integrable. Using our
main result, we check this for all two-dimensional examples, since
the complete integrability is relatively easy to work out in that
case.

One important property of the algebraically integrable operators
is that their eigenfunctions can be calculated explicitly (at
least, in principle). In section 5 we explain how to find the
Bloch eigenfunctions for a given integrable generalized Lam\'e
operator, in particular, for the elliptic Calogero--Moser problem.
As a result, we will see that the Bloch solutions are parametrized
by the points of an algebraic variety, which is a covering of a
product of elliptic curves (in a perfect agreement with the
situation in dimension one, due to Krichever \cite{K}). Let us
mention that for the elliptic Calogero--Moser problem (in the
$A_n$ case) the Bloch eigenfunctions were calculated by Felder and
Varchenko \cite{FV}. Our procedure is different and more general
(at cost of being less effective). We also explain how these Bloch
solutions can be used to construct the discrete spectrum
eigenstates for the Calogero--Moser problem.

In the last three sections of the paper we consider two particular
examples of the generalized Lam\'e operators in dimension $2$, for
which we make the formulas for the Bloch solutions very explicit.
The first example is
\begin{equation}\label{1}
L=-\Delta+2\wp(x)+2\wp(y)+4\wp(x-y)+4\wp(x+y)\,.
\end{equation}
This is a special case of the elliptic Calogero--Moser problem of
the $B_2$-type. Our second example is
\begin{gather}\label{2}
L=-\Delta +2\wp(x)+2(a^2+b^2)\wp(ax+by)+ 2(\overline{\mathstrut
a}^2+\overline{\mathstrut b}^2) \wp(\overline{\mathstrut
a}x+\overline {\mathstrut b}y)\,,\\\notag
a=\frac{-1+i\sqrt3\sin\alpha}{2}\,,\qquad \overline{\mathstrut
a}=\frac{-1-i\sqrt3\sin\alpha}{2}\,, \qquad
b=-\overline{\mathstrut b}=\frac{i\sqrt3\cos\alpha}2\,,
\end{gather}
where $\alpha$ is a complex parameter. Such operator was
considered by Hietarinta \cite{H} who showed that it admits a
commuting operator of order $3$. Its rational version
$\wp(z)=z^{-2}$ corresponds to a specific choice of the parameters
in the family of two-dimensional Schr\"odinger operators
introduced by Berest and Lutsenko \cite{BL} in connection with
Huygens' Principle (see section 4 of \cite{CFV} for details). Note
that the potential in \eqref{2} is real-valued (for real $x,y$) if
$\alpha$ is real and the period $\tau$ is pure imaginary.

It is more convenient to work with the following $3$-dimensional
version of \eqref{2}:
\begin{multline}\label{3}
L=-\d_1^2-\d_2^2-\d_3^2 +2(a_1^2+a_2^2)\wp(a_1x_1-a_2x_2)\\+
2(a_2^2+a_3^3)\wp(a_2x_2-a_3x_3)+2(a_3^2+a_1^2)\wp(a_3x_3-a_1x_1)\,,
\end{multline}
where $ a_1^2+a_2^2+a_3^2=0$, $\d_i=\d/\d x_i$. Then it is easy to
see that $L$ commutes with the operator $$ L_0=
a_1^{-1}\d_1+a_2^{-1}\d_2+a_3^{-1}\d_3 $$ and after restriction to
the plane $$ a_1^{-1}x_1 +a_2^{-1}x_2 +a_3^{-1}x_3 =0 $$ it
reduces to the operator \eqref{2} with proper $a,b$.

We calculate explicitly  the Bloch eigenfunctions of the operators
\eqref{1}, \eqref{3}. Notice a certain similarity between our
approach and the one used by Inozemtsev for $A_2$ case \cite{I}.
Let us also mention that in dimension two there is a nice theory
of Schr\"odinger operators which are finite-gap at a fixed energy
level, see \cite{DKN,VN}. It would be interesting to analyze our
results from that point of view.

In the last two sections we also calculate the Bloch solutions for
the discrete versions of \eqref{1}, \eqref{3}, which are given by
certain difference operators of Macdonald--Ruijsenaars type. This
raises a natural question about generalizing our results to the
difference setting. We hope to return to this problem in future.

\bigskip

\noindent{\bf Acknowledgments}. We are grateful to A.P.Veselov for
stimulating discussions and to Yu.Berest, S.Ruijsenaars and
K.Takemura for useful comments. P.E. thanks C. De Concini for a
discussion which was useful for the proof of Theorem
3.8. The work of O.C. was supported by
EPSRC. The work of A.O. was supported by Russian Foundation for
Basic Research (grant RFBR-01-01-00803).
The work of P.E. was partially supported by the NSF grant
DMS-9988796, and partly done for the Clay Mathematics Institute.

\section{Generalized Lam\'e operators}

Let $V=\mathbb C^n$ be a complex Euclidean space with the scalar
product denoted by $(\ ,\ )$, and $\mathcal A=\{\alpha\}$ be a
given finite set of affine-linear functions on $V$. Let us
consider a Schr\"odinger operator
\begin{equation}
\label{L} L=-\Delta+u(x)\,,\qquad \Delta=\partial^2/\partial
x_1^2+\dots +\partial^2/\partial x_n^2\,,
\end{equation}
with the elliptic potential $u$ of the following form:
\begin{equation}
\label{u} u(x)=\sum_{\alpha\in\mathcal
A}c_\alpha\wp(\alpha(x)|\tau)\,.
\end{equation}
Here $\wp(z|\tau)$ is the Weierstrass $\wp$-function with the
periods $1,\tau$, $\mathrm{Im}(\tau)>0$, and $c_\alpha\in\mathbb
C$ are the parameters which will be specified later (below we will
mostly suppress $\tau$, denoting $\wp(z|\tau)$ simply by
$\wp(z)$). Each of the functions $\alpha(x)$ in standard
coordinates on $V$ looks as $$\alpha(x)=a_0+a_1x_1+\dots
+a_nx_n\,,$$ so $\alpha$ is, effectively, a pair $(\alpha_0,a_0)$
where $a_0=\alpha(0)\in\mathbb C$ and
$\alpha_0=\mathrm{grad}\alpha=(a_1,\dots,a_n)$. We assume that
each $\alpha_0$ is non-isotropic (co)vector, i.e.
$(\alpha_0,\alpha_0)=a_1^2+\dots+a_n^2\ne 0$. Let us project
$\mathcal A$ onto $V^*$, $\alpha\mapsto\alpha_0$, denoting by
$\mathcal A_0$ the resulting set.

We want $u$ to be periodic, more precisely, we assume that the
lattice $\mathcal M\subset V^*$, generated over $\mathbb Z$ by the
set $\mathcal A_0$, has rank $\le n$. For simplicity, let us
assume that $\mathrm{rk}\mathcal M=n$ (the general case can be
reduced to that by passing to the factor space
$V/\mathrm{Ann}\mathcal M$). In that case the lattice $\mathcal
L+\tau\mathcal L$ is the period lattice for $u$, with $\mathcal
L:=\mathrm{Hom}(\mathcal M,\mathbb Z)\subset V$. Thus, $u$ may be
considered as a meromorphic function on a (compact) torus
$T=V/\mathcal L+\tau\mathcal L$, which is isomorphic to the
product of $n$ copies of the elliptic curve $\mathcal E=\mathbb
C/\mathbb Z+\tau\mathbb Z$. Singularities of $u(x)$ are the second
order poles along the following set of the hyperplanes: $$
Sing=\bigcup_{\alpha\in\mathcal A}\bigcup_{m,n\in\mathbb Z}
\pi_{\alpha}^{m,n}\,,\qquad\pi_\alpha^{m,n}:=\{x:\
\alpha(x)=m+n\tau\}\,.$$ We will assume that all hyperplanes
$\pi_\alpha^{m,n}$ are pairwise different; this can always be
achieved by rearranging the terms in \eqref{u}. This will imply
that $$u(x)\sim
c_\alpha(\alpha(x)-m-n\tau)^{-2}+O(1)\quad\text{near}\
\pi_\alpha^{m,n}\,.$$ Our next important assumption is that each
$c_\alpha$ in \eqref{u} has a form
$$c_\alpha=m_\alpha(m_\alpha+1)(\alpha_0,\alpha_0)\qquad\text{with
some}\quad m_\alpha\in\mathbb Z_{>0}\,.$$ Now we are going to put
some more restrictions on $u$ demanding its quasi-invariance in
the following sense.

\begin{defi} Let us say that the potential \eqref{u} as above is
{\bf quasi-invariant} if for any hyperplane
$\pi=\pi_\alpha^{m,n}\in{Sing}$ the (meromorphic) function
$u(x)-u(s_\pi x)$ is divisible by
$(\alpha(x)-m-n\tau)^{2m_\alpha+1}$, where $s_\pi$ denotes the
orthogonal reflection with respect to $\pi$.
\end{defi}

Here are two important examples of such potentials (more examples
will appear later).

\begin{example}\label{ex}
\label{cal} Consider the following potential $u(x)$ in $\mathbb
C^n$:
\begin{equation}
\label{A}
u=\sum_{i<j}2m(m+1)\wp(x_i-x_j)\,,\quad m\in\mathbb
Z_{>0}\,.
\end{equation}
It is invariant under any permutation of the coordinates
$x^1,\dots,x^n$, and also under translations $x\mapsto x+l$ for
$l\in\mathbb Z^n+\tau\mathbb Z^n$ ($\mathbb Z^n$ is the standard
integer lattice in $n$ dimensions). As a result, $u$ will be
symmetric with respect to any of the hyperplanes
$x_i-x_j=m+n\tau$. Thus, its Laurent expansion in the normal
direction will have no odd terms at all, hence $u$ is
quasi-invariant. The corresponding Schr\"odinger operator
\eqref{L} is the Hamiltonian of the quantum elliptic
Calogero--Moser problem. More generally, the elliptic
Calogero--Moser problems related to other root systems \cite{OP}
also lead, in the same way, to quasi-invariant potentials. In
particular, for the rank-one system $A_1$ we have the classical
Lam\'e operator \eqref{lame}
\end{example}

\begin{example}
\label{loc} This example is in dimension one. The potential $u$
has $N$ poles $x_1,\dots,x_N$ and looks as
\begin{equation}
\label{locus} u=\sum_{i=1}^N 2\wp(x-x_i)\,.
\end{equation}
 To ensure its
quasi-invariance, one has to impose the condition that $u$ has
zero derivative at each of its poles, more explicitly:
\begin{equation}\label{lloc}
\sum_{j\ne i}\wp'(x_i-x_j)=0\qquad\text{for all}\
i=1,\dots,N\,.
\end{equation}
This system of equations describes the so-called 'elliptic locus'
from \cite{AMM}, which has an intimate connection with the
classical elliptic Calogero--Moser system and the KdV hierarchy.
\end{example}

Let us call a Schr\"odinger operator $L$ with quasi-invariant
elliptic potential $u(x)$ a {\bf generalized  Lam\'e operator}. In
trigonometric and rational versions ($\wp(x)=\sin^{-2}x$ or
$x^{-2}$) such operators were considered in \cite{C98}, where
their eigenfunctions were effectively constructed. From the
results of \cite{C98} the so-called algebraic integrability of $L$
follows (see the paper \cite{CFV} for the rational case). We
recall the definition of the algebraic integrability in the next
section, following \cite{CV90,BEG}; let us just remark that in
dimension one this coincides with the class of algebro-geometric
operators which appear in the finite-gap theory, see \cite{KN,GW}.
This motivates the following

\medskip
\noindent{\bf Conjecture}.{\it The generalized Lam\'e operators
are all algebraically integrable.}
\medskip

As a particular case, this contains a conjecture of \cite{CV90}
about the algebraic integrability of the elliptic Calogero--Moser
problems. As we already mentioned in the introduction, for the
$A_n$-case \eqref{A} this has been proved by Braverman, Etingof
and Gaitsgory in \cite{BEG}. It is also known to be true in
dimension one, due to Gesztesy and Weikard \cite{GW1}. In the next
section we prove this conjecture under additional assumption of
the usual, Liouville integrability of $L$ (in a slightly stronger
sense). As a corollary, we will obtain the algebraic integrability
of the quantum Calogero--Moser problems for integer coupling
parameters.

\section{Monodromy and algebraic integrability}
\label{ai}

Let $L$ be a generalized Lam\'e operator as defined previously.
Recall that $L$ is {\bf completely integrable} if it is a member
of a commutative family of differential operators
$L_1=L,L_2,\dots,L_n$ which are algebraically independent. We
assume that the $L_i$'s have meromorphic coefficients and are
periodic with respect to the same lattice, which makes them
(singular) differential operators on the torus $T=\mathbb
C^n/\mathcal L+\tau\mathcal L$. The following proposition shows
that possible singularities of $L_i$ are contained in the singular
locus $Sing$ of the Schr\"odinger operator $L$.

\begin{prop} Let $L$ be a Schr\"odinger operator
regular in an open set $U$. Then any differential operator $M$ on
$U$ with meromorphic coefficients commuting with $L$ is regular in
$U$.
\end{prop}

To prove the proposition, we will need the following lemma.

\begin{lemma}
\label{ber} Let $S(x,p)$ be a meromorphic function on
$T^*U=U\times V^*$ which is polynomial in the momentum $p$, and
$\lbrace{ p^2,S(x,p)\rbrace}$ is a regular function. Then $S$ is a
regular function.
\end{lemma}
\begin{proof} Assume the contrary. The function $S$ is a finite
sum $\sum S_{\mathbf k}(x)p^\mathbf k$, where $p^\mathbf k$ are
monomials. Let $D\subset U$ be the divisor of poles of $S$. Take a
generic point $z_0$ of this divisor. Near this point $D$ is given
by an equation $f=0$, where $f$ is analytic at $z_0$ and
$df(z_0)\ne 0$. Let $S=\frac{Q}{f^k}(1+O(f))$ near $z_0$ ($Q$ is
regular at $z_0$, with $Q(z_0)\ne 0$). Then $$ \lbrace
p^2,S\rbrace= -2k\sum p_i\frac{\partial f}{\partial
x_i}f^{-k-1}Q+O(f^{-k})\,.$$ But $\sum p_i\frac{\partial
f}{\partial x_i}(z_0,p)\ne 0$ for generic $p$. Thus, $\lbrace
p^2,S\rbrace$ is singular, which is a contradiction.
\end{proof}

Now we prove the proposition. Suppose $[M,L]=0$. Assume $M$ is not
regular. Then we can write $M$ as $M'+M''$, where $M''$ is
regular, and $M'$ has a singular highest symbol. Then
$[L,M']=-[L,M'']$ is regular. Let $S(x,p)$ be the symbol of $M'$.
We have $\lbrace{ p^2,S\rbrace}$ is regular (since if it is
nonzero, it is the symbol of $[L,M']$). Then the lemma implies
that $S$ is regular. This contradiction proves the proposition.

Notice also that if $S(x,p)$ is the highest symbol of $L_i$, then
from $[L,L_i]=0$ it follows that $\lbrace p^2,S\rbrace =0$. Thus,
by Lemma \ref{ber}, $S$ must be regular everywhere. Hence, each
$L_i$ must have constant highest symbol, i.e. $L_i=s_i(\partial)+
\dots$ for some polynomial $s_i$.

\begin{defi} Let us say that a Schr\"odinger operator
$L=-\Delta+u$ in $V=\mathbb C^n$ is {\bf strongly integrable} if
the commuting operators $L_1=L,\dots,L_n$ have algebraically
independent homogeneous constant highest symbols $s_1,\dots,s_n$
and if $\mathbb C[V]$ is finitely generated as a module over the
ring generated by $s_1,\dots,s_n$ or, equivalently, if the system
$s_1(\xi)=0,\dots, s_n(\xi)=0$ has the unique solution $\xi=0$.
\end{defi}

Now let $L$ be a strongly integrable generalised Lam\'e operator,
so we have the operators $L_1=L,\dots, L_n$ with meromorphic
coefficients on the torus $T=\mathbb C^n/\mathcal L+\tau\mathcal
L$, and $L_i=s_i(\partial)+\dots$. First of all, $\mathbb C[V]$ is
locally free as a module over $\mathbb C[s_1,\dots, s_n]$. This
follows from the fact, due to Serre \cite{S}, that if $f:X\to Y$
is a finite map of smooth affine varieties of the same dimension,
then $\mathcal O(X)$ is a locally free $\mathcal O(Y)$-module.
Further, since $s_i$ are homogeneous, this module is graded,
hence, it must be free. Denote by $N$ the rank of this free
module.

Consider now the eigenvalue problem
\begin{equation}
\label{pro} L_i\psi=\lambda_i\psi\,,\quad i=1,\dots, n\,,
\end{equation}
with $\lambda=(\lambda_1,\dots, \lambda_n)\in\mathbb C^n$. Then
the space of solutions of this system in any simply connected
domain in $T\setminus Sing$ is N-dimensional. So we have a
holonomic system of rank $N$ on $T$ with singularities along a
finite union of hypertori $Sing\subset T$.

\begin{theorem} The holonomic system \eqref{pro} has regular singularities.
\end{theorem}

\begin{proof}
Regularity of singularities is a codimension one condition. Thus,
it is sufficient to restrict our attention to a neighborhood $U$
of a point which lies on exactly one subtorus from the pole
divisor. Let us assume that this point is $0$, and the subtorus is
locally defined by the equation $x_1=0$.

\begin{lemma}\label{bound}
For any $r\ge 0$, the singular part of the differential operator
$(x_1)^rL_i$ has the order at most $d_i-r-1$, where $d_i$ is the
degree of $s_i$.
\end{lemma}

\begin{proof} Let us introduce new variables $y_j=tx_j$, and
write our operators with respect to them. Then
$L=t^{-2}L^{(0)}+O(1)$, $t\to 0$, where
$L^{(0)}=-\Delta+\frac{m_\alpha(m_\alpha+1)}{x_1^2}$ with
$m_\alpha\in\mathbb Z_+$ corresponding to the chosen hyperplane
$x_1=0$. Let $L_i^{(0)}$ be the coefficient at the lowest power of
$t$ in the expression for $L_i$. Then $L_i^{(0)}$ is homogeneous,
of degree $\le -d_i$ in $t$ (since the symbol of $L_i$ has this
degree). On the other hand, the order of this operator is at most
$d_i$. Thus, if this degree of $L_i^{(0)}$ is less than $d_i$,
then its symbol would have to be singular. But $L_i^{(0)}$
commutes with $L^{(0)}$, so this contradicts Lemma \ref{ber}.
Thus, the degree is exactly $d_i$, which implies the statement.
\end{proof}

\begin{lemma}\label{regu}
Consider a system of differential equations $df/dz=A(z)f(z)$ with
holomorphic coefficients on a punctured disk $0<|z|<\epsilon$ ($f$
is a vector function, $A$ is a matrix function). Assume that
$A(z)$ is meromorphic at $z=0$ and that there exists an
integer-valued function $D(i)$ such that the order of $a_{ij}(z)$
at $z=0$ is at least $D(j)-D(i)-1$. Then the system has a regular
singularity at $z=0$.
\end{lemma}

Indeed, let us make a change of variable $g_i=z^{D(i)}f_i$. This
will change the matrix $A$ into a new matrix $\tilde A$ which
obviously has at most simple pole at $z=0$.

Now let us prove the theorem. First, let us reduce \eqref{pro} to
a first order holonomic system in a standard way. Choose a
collection of $N$ homogeneous polynomials $q_1,...,q_N$ in
$\partial_i$, which form a basis in $\mathbb
C[\partial_1,...,\partial_n]$ as a free module over $\mathbb
C[s_1(\partial),...,s_n(\partial)]$. Let $\psi$ be a solution the
eigenvalue problem \eqref{pro}, and consider the functions
$$q_1(\partial)\psi,...,q_N(\partial)\psi\,.$$ Then, for any
polynomial $q$ of $\partial_i$, the function $q\psi$ can be
expressed via $q_i\psi$, with the coefficients depending on $x$
and $\lambda$. Indeed, assume that $q$ is homogeneous of degree
$d$. We know that $q$ can be uniquely represented in the form
$\sum g_jq_j$, where $g_j\in \mathbb C[s_1,...,s_n]$. Thus,
$q\psi-\sum q_jg_j(L_1,..,L_n)\psi$ is expressed through
differential polynomials of $\psi$ of degree smaller than $d$. But
$g_j(L_1,...,L_n)\psi=g_j(\lambda_1,...,\lambda_n)\psi$, so,
eventually, by induction we get a desired representation of
$q\psi$.

Now observe that by Lemma \ref{bound}, the coefficient of
$q_j\psi$ in the expansion of $q\psi$ has a pole in $x_1$ of order
at most $\deg(q)-\deg(q_j)$. Thus, we have a holonomic system of
matrix partial differential equations $$\partial_kq_i\psi=\sum
a_{ij}^{(k)}(x)q_j\psi\,\quad k=1,\dots, n\,,$$ which is
equivalent to the system \eqref{pro}. Each of the matrices
$a^{(k)}_{ij}$ ($k=1,\dots, n$) satisfies the conditions of Lemma
\ref{regu} with respect to $z=x_1$, for $D(i)=-\deg(q_i)$. In
particular, using the lemma, we conclude from the first equation
$\partial_1q_i\psi=\sum a_{ij}^{(1)}\psi$ that all the solutions
have at most power growth in $x_1$ when approaching $x_1=0$, hence
the system \eqref{pro} has a regular singularity at $x_1=0$.
\end{proof}

Now take any point $x_0\in T\setminus Sing$ and consider the
monodromy of the system \eqref{pro} at point $x_0$, this gives an
$N$-dimensional representation of the fundamental group
$\pi_1(T\setminus Sing)$.

\begin{prop} The monodromy group of the system \eqref{pro} is
commutative for any $\lambda$.
\end{prop}
\begin{proof}
Take a hyperplane $\pi=\pi_\alpha^{m,n}\in{Sing}$ and let $P$ be a
generic point of $\pi$. Changing the coordinates if necessary, we
may assume that $P=(0,\dots,0)$ is the origin and $\pi$ is given
by equation $x_1=0$. From the regularity of singularities it
follows that there exist $\gamma_1,\dots, \gamma_r\in\mathbb C$
such that any solution $\psi$ of the system \eqref{pro} near
$P\in\pi$ has a convergent series expansion in the subspace
$$\bigoplus_{j=1}^r x_1^{\gamma_j}\mathbb
C[[x_1,\dots,x_n]][\log(x_1)]\,.$$ Substituting such a series into
the first equation $L\psi=\lambda_1\psi$, one arrives at certain
recurrence relations from which follows that (1) we have two
'leading exponents' $\gamma_1=-m_\alpha$, $\gamma_2=m_\alpha+1$,
both integer; (2) there will be no $\log(x_1)$ involved. The
latter fact is due to the quasi-invariance of $u$, see \cite{DG}
for the one-dimensional case and \cite{CFV}, section 2 for a
discussion in the multivariable setting.

As a result, we see that all the solutions are single-valued near
$\pi$, so the local monodromy corresponding to a loop around $\pi$
is trivial. Consequently, the global monodromy group will be a
homomorphic image of the commutative group $\pi_1(T)$.
\end{proof}

\begin{cor} The differential Galois group of the system
\eqref{pro} is commutative for any $\lambda$.
\end{cor}
\begin{proof}
It is known (see \cite{De}) that for a regular holonomic system on
a smooth projective variety the differential Galois group
coincides with the Zariski closure of the monodromy group. In our
situation the proof is simple: first, we know that all solutions
of the system \eqref{pro} are meromorphic in $\mathbb C^n$.
Indeed, this is true outside the set of codimension $2$, by the
proposition above, hence it holds everywhere by Hartogs' theorem.
Now let $F$ be the field of 'elliptic functions', i.e. meromorphic
functions on the torus $T$, and $L$ is the solution field of
\eqref{pro} (on some simply connected domain). The monodromy gives
the homomorphism $$\pi_1(T)\mapsto {Dgal}(L:F)$$ of the
fundamental group of the torus $T$ into the differential Galois
group. Let $G$ be the Zariski closure of the image of this
homomorphism inside ${Dgal}(L:F)$. By the main theorem of the
differential Galois theory \cite{Kap}, to prove that
$G={Dgal}(L:F)$ it suffices to show that in the solution field any
$G$-invariant function is actulaly ${Dgal}$-inavriant, i.e.
belongs to $F$. But any meromorphic $G$-invariant function is
$\pi_1(T)$-invariant, hence elliptic. Thus, it is
${Dgal}$-invariant by definition.

This proves that ${Dgal}(L:F)=G$, and the latter is obviously
commutative.
\end{proof}

\medskip

Now let us remark (see \cite{BEG}) that a quantum completely
integrable system (QCIS) on a smooth algebraic variety $X$ of
dimension $n$ naturally defines an embedding $$\theta: \mathcal
O(\Lambda)\mapsto D(X)\,,$$ where $\Lambda\simeq \mathbb A^n$ is
an affine space and $D(X)$ denotes the ring of differential
operators on $X$. More generally, $\Lambda$ can be any affine
variety with $\dim(\Lambda)=n$. Then we have an analogous
eigenvalue problem $\theta(g)\psi=g(\lambda)\psi$, $\forall
g\in\mathcal O(\Lambda)$. The dimension of the local solution
space of this system at generic point of $X$ is called the {\bf
rank} of a QCIS. Recall further, that a QCIS $S=(\Lambda, \theta)$
is {\bf algebraically integrable} if it is dominated by another
QCIS $S'=(\Lambda',\theta')$ of {\bf rank one} ($S$ is dominated
by $S'$ if there is a map of algebras $h: \mathcal
O(\Lambda)\mapsto \mathcal O(\Lambda')$ such that
$\theta=\theta'\circ h$). In our situation, this implies that
apart from the operators $L_1,\dots, L_n$, we have additional
commuting operators which are not algebraic combinations of $L_i$
(though, of course, they are algebraically dependent with $L_i$).
In dimension $n=1$ this is equivalent to saying that $L$ is a
member of a maximal commutative ring in $D(X)$ of rank one; this
is known to coincide with the class of algebro-geometric
operators.

\begin{theorem}\label{aii} Any generalised Lam\'e operator $L$ which is strongly
integrable is algebraically integrable.
\end{theorem}
This follows immediately from the result above and the criterion
from \cite{BEG}. Moreover, according to \cite{BEG}, for generic
$\lambda$ the solution space is generated by the quasiperiodic
solutions:
\begin{cor} There exist meromorphic $1$-forms $\omega_j$
on the torus $T$, with first order poles and depending
analytically on $\lambda$, such that the functions $\psi_j=e^{\int
\omega_j}$ give a basis of the solution space of \eqref{pro} for
generic $\lambda$.
\end{cor}
Each of these functions will be double-Bloch, in terminology of
\cite{KZ}:
\begin{equation}\label{qua}
\psi_j(x+l)=e^{2\pi i\langle a_j, l\rangle}\psi_j(x)\,,\qquad
\psi_j(x+\tau l)=e^{2\pi i\langle b_j, l\rangle}\psi_j(x)\,,
\end{equation}
for appropriate $a_j, b_j\in V^*$ and for all $l\in\mathcal L$.
Namely, $$\langle a_j, l\rangle =\frac{1}{2\pi i}
{\int_z^{z+l}\omega_j}\quad\text{and}\ \langle b_j,
l\rangle=\frac{1}{2\pi i}{\int_z^{z+\tau l}\omega_i}$$ (since
$\psi_j$ has no branching along ${Sing}$, these are well defined
modulo $\mathbb Z$). In Section \ref{bloch} below we explain how
one can calculate these double-Bloch solutions for the generalized
Lam\'e operators.

\begin{remark}
Our argument applies to a more general situation when one has
commuting operators $L_1,\dots,L_n$ on an abelian variety, such
that the system $L_i\psi=\lambda_i\psi$, $i=1,\dots ,n$, is
regular holonomic. Then the triviality of its local monodromy
around singularities implies that this system is algebraically
integrable.
\end{remark}

\section{Examples}

\subsection{One-dimensional case}

Let $L$ be a one-dimensional Schr\"odinger operator
$L=-\frac{d^2}{dx^2}+u(x)$, $x\in\mathbb C$. Consider the
eigenvalue problem $L\psi=\lambda\psi$. Now, assuming that $u$ is
meromorphic and has a pole at $x=0$, let us consider the local
monodromy of the solutions of this second order differential
equation. Suppose that this monodromy is trivial {\it for all}
$\lambda$, in other words, all the solutions are meromorphic at
$x=0$. Then, using the classical Frobenius analysis, one can see
that $u$ must have a pole of the second order, with no residue:
$u=m(m+1)x^{-2}+O(1)$, with integer $m$. Furthermore, in the
series for $u$ at $x=0$ there must be no terms of order $2j-1$ for
all $j=1,\dots,m$ (see \cite{DG}). This is exactly the
quasi-invariance of $u$ with respect to the symmetry $x\to -x$.

Now let $u$ be elliptic, with periods $1,\tau$. Let us demand all
the solutions of the equation $L\psi=\lambda\psi$ to be
meromorphic in $\mathbb C$ (such $u$ are called Picard potentials
in \cite{GW,GW1}). By the discussion above, this is equivalent to
the quasi-invariance of $u$ at each pole. More explicitly,
$$u(x)=\sum_{i=1}^Mm_i(m_i+1)\wp(x-x_i)\,,\quad\text{with}\quad
m_i\in\mathbb Z_+\,,$$ and the poles $x_i$ must satisfy the
following system of equations, which generalizes \eqref{lloc}:
\begin{equation}\label{loc1}
\sum_{j\ne i}m_j(m_j+1)\wp^{(2s-1)}(x_i-x_j)=0
\end{equation}
for all $i=1,\dots, M$ and $s=1,\dots, m_i$.

Such $L$ automatically defines a $QCIS$ (of rank 2), and the
regularity of singularities is obvious. Applying the results of
the previous section, we conclude that $L$ is algebraically
integrable. Algebraic integrability of $L$ implies the existence
of a differential operator $P$, commuting with $L$. Since $P$ is
not a polynomial of $L$, we may assume that it is of odd order.
Thus, $L$ must be algebro-geometric, according to the
Burchnall--Chaundy--Krichever theory, see e.g.\cite{SW}. Thus, our
result in this case is equivalent to the main result of
\cite{GW1}. Notice that our approach easily extends to the case of
operators of any order (cf. the remark at the end of \cite{GW1}).

\subsection{Quantum Calogero--Moser system}

Let $V$ be a complex Euclidean space, $\dim V=n$. Let
$R=\{\alpha\}$ be a reduced irreducible root system in $V^*$, $W$
be the corresponding Weyl group, and the parameters $m_\alpha$ be
chosen in a $W$-invariant way. The corresponding Calogero--Moser
operator \cite{OP} looks as
\begin{equation}\label{cmroot}
L=-\Delta+\sum_{\alpha\in
R_+}m_\alpha(m_\alpha+1)(\alpha,\alpha)\wp(\langle\alpha,x\rangle)\,.
\end{equation}
Let $\mathbb C[V]^W$ be the ring of $W$-invariant polynomials. By
the Chevalley theorem, it is freely generated by $n$ elements
$p_1,\dots,p_n$, and $\mathbb C[V]$ is a free module over $\mathbb
C[V]^W$, of rank $|W|$. The following result has been proved in
\cite{Ch} for any $W$-invariant $m_\alpha\in\mathbb C$.

\begin{theorem}[Cherednik]\label{ch} For each homogenious $p\in\mathbb C[V]^W$
there exists a differential operator $L_p$ with the highest symbol
$p$, commuting with $L$: $[L,L_p]=0$. The family $\{L_p\}$ is
commutative.
\end{theorem}

For integer $m_\alpha$ the potential $u$ in \eqref{cmroot} is
quasi-invariant (cf. Example \ref{ex}). Altogether this proves the
conjecture of \cite{CV90}:
\begin{cor} The Calogero--Moser operator \eqref{cmroot} is
algebraically integrable for integer $m_\alpha$.
\end{cor}

Let us mention that for the classical root systems $R=A\dots D$
the complete set of commuting operators $L_1,\dots, L_n$ for
\eqref{cmroot} was explicitly found by Ochima and Sekiguchi
\cite{OS}. Their results cover the $BC_n$ case, too. In that case
we have $5$ parameters $m,g_0,g_1,g_2,g_3$ and the Inozemtsev
operator
\begin{equation}\label{bcn}
L=-\Delta+2m(m+1)\sum_{i<j}(\wp(x_i-x_j)+\wp(x_i+x_j))+
\sum_{i=1}^n\sum_{s=0}^3 g_s(g_s+1)\wp(x_i+\omega_s)\,,
\end{equation}
with $\omega_s (s=0\dots 3)$ denoting the half-periods
$0,1/2,\tau/2,(1+\tau)/2$. The Weyl group $W$ for this case is
generated by the permutations of $x_i$ and sign flips, and Theorem
\ref{ch} still holds true, due to \cite{OS}. Applying Theorem
\ref{aii} we obtain

\begin{cor} The Inozemtsev operator \eqref{bcn} is
algebraically integrable for any integer $m, g_0, g_1, g_2, g_3$.
\end{cor}

\subsection{Deformed root systems}

Other known examples of the generalized Lam\'e operators in
dimension $>1$ are related to deformed root systems, which
appeared in \cite{CFV}. Below we describe the set of linear
functionals $\mathcal A=\{\alpha\}$ and the corresponding
multiplicities $m_\alpha$.

\medskip
\noindent (1) $A_{n,1}(m)$ system \cite{CFV}.

\noindent It consists of the following covectors in ${\mathbb
C}^{n+1}$: $$
%A_{n-1,1}(m)=
\left\{
\begin{array}{lll}
x_i - x_j, &  1\le i<j\le n, & \text{with multiplicity } \langle
m\rangle\,,\\ x_i - \sqrt{m}x_{n+1}, &  i=1,\ldots ,n & \text{with
multiplicity }  1\,.
\end{array}
\right. $$ Here $m$ is an integer parameter, and $\langle m
\rangle$ stands for $\langle m \rangle=\max\{m,-1-m\}$.

\medskip
\noindent (2) $C_{n,1}(m,l)$ system \cite{CFV}.

\noindent It consists of the following covectors in $\mathbb
C^{n+1}$: $$ \left\{
\begin{array}{lll}
x_i\pm x_j, &  1\le i<j\le n, &  \text{with multiplicity } \langle
k\rangle\,,\\ 2x_i, &  i=1,\ldots ,n  & \text{with multiplicity }
\langle m\rangle\,,\\ x_i\pm \sqrt{k}x_{n+1}, & i=1,\ldots ,n &
\text{with multiplicity } 1\,,\\ 2\sqrt{k}x_{n+1} && \text{with
multiplicity } \langle l\rangle\,.\\
\end{array}
\right. $$ Here $k,l,m$ are integer parameters related as $k
=\frac{2m+1}{2l+1}$, and $\langle k\rangle, \langle l\rangle,
\langle m\rangle$ have the same meaning as in $A_{n,1}(m)$ case.
In two-dimensional case $n=1$ the first group of roots is absent
and there is no restriction for $k$ to be integer.

\medskip
\noindent (3) Here is a $BC_n$-type generalization of the previous
example. Let $\omega_s (s=0\dots3)$ denote the half-periods, as in
$BC_n$ case \eqref{bcn}. The set of linear functionals
$\alpha\in\mathcal A$ and the corresponding multiplicities look as
follows: $$ \left\{
\begin{array}{lll}
x_i\pm x_j, &  1\le i<j\le n, &  \text{with multiplicity } \langle
k\rangle\,,\\ x_i+\omega_s, &  i=1,\ldots ,n  & \text{with
multiplicity } \langle m_s\rangle\,,(s=0\dots 3)\,,\\ x_i\pm
\sqrt{k}x_{n+1}, & i=1,\ldots ,n & \text{with multiplicity }
1\,,\\ \sqrt{k}x_{n+1}+\omega_s && \text{with multiplicity }
\langle l_s\rangle\,, (s=0\dots 3)\,.\\
\end{array}
\right. $$ Here $k,l_s,m_s$ are nine integer parameters related
through $k =\frac{2m_s+1}{2l_s+1}$ for all $s=0\dots 3$. The
previous case corresponds to $m_s\equiv m$ and $l_s\equiv l$.
Again, in case $n=1$ the first group of roots is absent and $k$
may not be integer.

\medskip
\noindent(4) Hietarinta operator \cite{H}.

\noindent In this case we have three covectors in $\mathbb C^3$,
$$\alpha=a_1x_1-a_2x_2\,,\beta=a_2x_2-a_3x_3\,,\gamma=a_1x_1-a_3x_3\,,$$
with $m_\alpha=m_\beta=m_\gamma=1$. Here $a_i$ are arbitrary
complex parameters such that $a_1^2+a_2^2+a_3^2=0$. Notice that
the system is essentially two-dimensional since
$\alpha+\beta+\gamma=0$.

\medskip
\noindent(5) $A_{n-1,2}(m)$ system \cite{CV1}.

\noindent It consists of the following covectors in $\mathbb
C^{n+2}$: $$\left\{
\begin{array}{lll}
x_i - x_j, &  1\le i<j\le n, & \text{with multiplicity } m\,,\\
x_i - \sqrt{m}x_{n+1}, &  i=1,\ldots ,n  & \text{with multiplicity
} 1\,,\\ x_i - \sqrt{-1-m}x_{n+2}, & i=1,\ldots ,n & \text{with
multiplicity }  1\,,\\ \sqrt{m}x_{n+1}-\sqrt{-1-m}x_{n+2} &  &
\text{with multiplicity }  1\,.
\end{array}
\right. $$ Notice that for $m=1$ this system coincides with the
system $A_{n,1}(-2)$ above.

\medskip
In all these cases a direct check shows that the corresponding
potential $u$ is quasi-invariant. Notice that in all cases $u$ is
symmetric with respect to $s_\alpha$ as soon as $m_\alpha>1$.
Thus, one have to check the quasi-invariance only for those
$\alpha$ where $m_\alpha=1$.

We believe that all these operators are algebraically integrable.
In the next section we check this for all two-dimensional
examples.

\subsection{Two-dimensional case}

Let us check that all known two-dimensional generalized Lam\'e
operators are algebraically integrable. Apart from the root
systems $A_2, BC_2$ and $G_2$ considered previously, we have three
deformed cases, namely, the $A_{1,1}(m)$ case, the Hietarinta
operator and the deformed $BC_2$ case.

First, let us consider the $A_{1,1}(m)$ case:
\begin{multline}\label{21}
L=-\d_1^2-\d_2^2-\d_3^2+2m(m+1)\wp(x_1-x_2)+\\2(m+1)\wp(x_1-\sqrt
mx_3)+2(m+1)\wp(x_2-\sqrt mx_3)\,.
\end{multline}
In this case we can use the result of \cite{Xo} where the complete
integrability  of \eqref{21} was established. First, $L$ obviously
commutes with $L_0=\d_1+\d_2+\frac{1}{\sqrt m}\d_3$.
\begin{prop}[\cite{Xo}] For any $m$ there exists a third order operator
$L_2$ commuting with $L_0,L$, and its highest symbol is
$\d_1^3+\d_2^3+\sqrt m\d_3^3$.
\end{prop}
It is easy to check that as soon as $m\ne -1$, the highest symbols
of $L_0,L,L_2$ will satisfy the requirements of the strong
integrabilty (notice that for $m=-1$ the operator $L$ is trivial).
As a result, for integer $m$ we obtain the algebraic integrability
of the operator \eqref{21}. Note that for the special case $m=2$
the algebraic integrability of \eqref{21} was demonstrated in
\cite{Xo1} by presenting an explicit extra operator commuting with
$L_0, L, L_2$.

Now let us consider the Hietarinta operator
\begin{multline}\label{22}
L=-\d_1^2-\d_2^2-\d_3^2+2(a_1^2+a_2^2)\wp(a_1x_1-a_2x_2)+\\
2(a_2^2+a_3^3)\wp(a_2x_2-a_3x_3)+2(a_3^2+a_1^2)\wp(a_3x_3-a_1x_1)\,,
\end{multline}
where $a_1^2+a_2^2+a_3^2=0$. Such $L$ commutes with
$L_0=a_1^{-1}\d_1+a_2^{-1}\d_2+a_3^{-1}\d_3$, and we need one more
operator for the complete integrability (we assume that all $a_i$
are nonzero, otherwise $L$ is reducible). Such operator was found
in \cite{H}.
\begin{prop}[\cite{H}] There exists a third order operator $L_2$ commuting
with $L_0,L$ above, and its highest symbol is
$a_1\d_1^3+a_2\d_2^3+a_3\d_3^3$.
\end{prop}
If we denote the highest symbols of these three operators as
$s_1,s_2,s_3$, then it is easy to check that the system
$s_1(\xi)=s_2(\xi)=s_3(\xi)=0$ has the only solution $\xi=0$; the
only exception is the case when $a_1^3=a_2^3=a_3^3$. As a result,
we conclude that for all values of the parameters (apart from the
case $a_1^3=a_2^3=a_3^3$) the Hietarinta operator is strongly
integrable, thus, it is algebraically integrable. Note that this
also follows from \cite{H} where one more operator commuting with
$L_0, L, L_2$ was found.

\medskip
\begin{remark} In the case $a_1^3=a_2^3=a_3^3$ the operator \eqref{22} is still
algebraically integrable, although it is no longer strongly
integrable.
\end{remark}

\medskip
Finally, let us consider the deformed $BC_2$ case. The
Schr\"odinger operator $L$ has the following form:
\begin{equation}\label{23}
L=-\d_x^2-\d_y^2+U(x,y)\,,
\end{equation}
where $U=2(k+1)(\wp(x+\sqrt{k}y)+\wp(x-\sqrt{k}y))+v(x)+w(y)$ and
$v,w$ are given by the expressions
\begin{equation}\label{23a}
v=\sum_s m_s(m_s+1)\wp(x+\omega_s)\,,\quad w=k\sum_s
l_s(l_s+1)\wp(\sqrt{k}y+ \omega_s)\,.
\end{equation}
Here $\o_0,\o_1,\o_2,\o_3$ are the half-periods and $m_s,l_s$ and
$k$ are nine parameters such that $k=(2m_s+1)/(2l_s+1)$ for all
$s=0,1,2,3$ (thus, effectively, $L$ contains five independent
parameters).

\begin{prop} For any values of the parameters $m_s,l_s,k$ such that
$k=\frac{2m_s+1}{2l_s+1}$ the following operator commutes with
$L$:
\begin{align*}
M=&-\partial_x^4-k\partial_y^4+2(U-w)\partial_x^2-
4\sqrt{k}(k+1)(\wp(x+\sqrt{k}y)-\wp(x-\sqrt{k}y))\partial_x\partial_y+\\
&2k(U-v)\partial_y^2 +(2v'+2(k+1)(2-k)(\wp'(
x+\sqrt{k}y)+\wp'(x-\sqrt{k}y)))\partial_x+\\
&(2k^2w'+2\sqrt{k}(k+1)(2k-1)(\wp'(x-\sqrt{k}y)-\wp'(x+\sqrt{k}y)))
\partial_y-\\
&(k+1)^3\wp(x+\sqrt{k}y)\wp(x-\sqrt{k}y)+8(k^3+1)
(\wp^2(x+\sqrt{k}y)+\wp^2(x-\sqrt{k}y))-\\&4(k+1)(v+kw)
(\wp(x+\sqrt{k}y)+\wp(x-\sqrt{k}y))+v''-v^2+ k(w''-w^2)\,.
\end{align*}
Here $v',w'$ and so on are the derivatives with respect to the
corresponding variable.
\end{prop}
Now since the system $\xi_1^2+\xi_2^2=\xi_1^4+k\xi_2^4=0$ does not
have nontrivial solutions as soon as $k\ne-1$, we conclude that
$L$ is strongly integrable (for $k=-1$ this is also true because
$L$ is reducible in that case). Hence, for any integer values of
the parameters $l_s,m_s$ the operator \eqref{23}--\eqref{23a} is
algebraically integrable.

\section{Bloch solutions}
\label{bloch}

Let $L$ be a generalized Lam\'e operator which is strongly
integrable, thus algebraically integrable. We know already that
for generic $\lambda$ the solution space of \eqref{pro} is spanned
by the meromorphic double-Bloch solutions. Now we are going to
explain how one can, in principle, calculate them.

Let $\mathcal W$ denote the following linear subspace in the space
of meromorphic functions on $\mathbb C^n$. First, its elements are
holomorphic everywhere apart from the singular locus $Sing=\bigcup
\pi_\alpha^{m,n}$ of the operator $L$, where they may have poles,
of order $\le m_\alpha$ along $\pi_\alpha^{m,n}$. Next, take any
hyperplane $\pi=\pi_\alpha^{m,n}$ with $s_\pi$ denoting the
orthogonal reflection with respect to $\pi$. Then any function
$\varphi\in \mathcal W$ must have the following property:
\begin{equation}\label{qq}
\varphi(x)-(-1)^{m_\alpha}\varphi(s_\pi x)\ \text{ is divisible
by}\quad (\alpha(x)-m-n\tau)^{m_\alpha+1}\,.
\end{equation}

\begin{prop}\label{f} The subspace $\mathcal W$ is stable under the action of
$L$: $L(\mathcal W)\subseteq\mathcal W$. Furthermore, any
meromorphic eigenfunction $\varphi$ of $L$ must belong to
$\mathcal W$. The same is true for any of the commuting operators
$L_2,\dots,L_n$: $L_i(\mathcal W)\subseteq\mathcal W$.
\end{prop}

\begin{proof} Consider any
hyperplane $\pi\in Sing$ and adjust the coordinates in such a way
that $\pi$ is given by equation $x_n=0$. Take a generic point in
$\pi$ and expand $\varphi$ in Laurent series in normal direction
to $\pi$, i.e. $\varphi=\sum_{j\in\mathbb Z}a_j(x_n)^j$,
$a_j=a_j(x_1,\dots,x_{n-1})$. Now put $$M:=\{-m_\alpha+2\mathbb
Z_{\ge 0}\}\cup\{m_\alpha+1+2\mathbb Z_{\ge 0}\}$$ and split the
series into two parts, with $j\in M$ and with $j\in\mathbb Z
\backslash M$: $\varphi=\varphi_1+\varphi_2$. First claim is that
an application of $L$ to $\varphi_1$ will produce a series of a
similar kind. This follows directly from the quasi-invariance of
$u$, proving the first part of the proposition. On the other hand,
if $a_{j}(x_n)^j$ is the first nonzero term in $\varphi_2$, then
applying $L$ to $\psi_2$ will give a series starting from
$(x_n)^{j-2}$, which would contradict the equation
$L\varphi=\lambda\varphi$, thus proving the second claim.

In a similar way, if $L'L=LL'$ for some other operator $L'$, then
$L'(L)^r=(L)^rL'$ for any $r\ge 1$. Thus, if $\mathcal
W':=L'(\mathcal W)$ then $L^r(\mathcal W')\subseteq \mathcal W'$
for any $r$. Now suppose we could find a function $\varphi$ in
$\mathcal W'$ which is not in $\mathcal W$. Then, consider a
series expansion of $\varphi$ in the direction, normal to $\pi$,
and split it into $\varphi=\varphi_1+\varphi_2$ as above. If
$\varphi_2\neq 0$, then $L^r\varphi$ would have a pole of an
arbitrarily high order along $\pi$ (as $r$ increases), which is
impossible for an element in $\mathcal W'$ (since $\mathcal
W'=L'\mathcal W$ and $L'$ has meromorphic coefficients). This
would contradict the inclusion $L^r(\mathcal W')\subseteq \mathcal
W'$. Thus, $\mathcal W'=\mathcal W$.
\end{proof}

Now let $\psi$ be a double-Bloch solution of \eqref{pro}, so for
appropriate $a,b\in V^*$ and for any $l\in\mathcal L$ we have:
\begin{equation}
\label{char} \psi(x+l)=\psi(x)e^{2\pi i\langle a,
l\rangle}\,,\quad \psi(x+\tau l)=\psi(x)e^{2\pi i\langle b,
l\rangle}\,.
\end{equation}
We know that $\psi$ is meromorphic in $\mathbb C^n$ with possible
poles along the hyperplanes $\pi_\alpha^{m,n}$, of order
$m_\alpha$. In our discussion below we restrict ourselves to the
case when all the linear functions $\alpha\in\mathcal A$ have zero
constant term, i.e. $\alpha=\alpha_0\in V^*$ for all $\alpha$, so
$\alpha(x)=\langle\alpha,x\rangle$. Everything extends to the
general case with obvious modifications.

As a result, we see that $\psi$ can be presented in the form
\begin{equation}
\label{gans} \psi=\Phi/\delta\,,\qquad
\delta(x)=\prod_{\alpha\in\mathcal
A}\theta(\langle\alpha,x\rangle)^{m_\alpha}\,,
\end{equation}
for some holomorphic $\Phi(x)$. Here $\theta=\theta_1$ is the
classical (odd) Jacobi theta function,
\begin{equation}\label{theta1}
\theta(z)= \sum_{n\in{\mathbb Z}}\exp(\pi i (n+1/2)^2\tau+2\pi
i(n+1/2)(z+1/2))\,.
\end{equation}
Recall that $\theta(z)$ has the following translation properties
in $z$: $$\theta(z+1)=-\theta(z)\,,\qquad\theta(z+\tau)=-e^{-2\pi
iz-\pi i\tau}\theta(z)\,.$$ This determines the translation
properties of $\delta$. To write them down, it is convenient to
introduce the following linear map $\Omega: V\to V^*$ which is
defined as
\begin{equation}\label{Om}
\Omega: x\mapsto\sum_{\alpha\in\mathcal
A}m_\alpha\langle\alpha,x\rangle\alpha\,.
\end{equation}
Note that $\Omega$ maps the lattice $\mathcal L$ to a sublattice
in $\mathcal M=\mathrm{Hom}(\mathcal L,\mathbb Z)$. We also need a
covector $\varrho=\frac12\sum_{\alpha\in\mathcal
A}m_\alpha\alpha$.

Under these notations, we have the following translation formulas
for any $l\in\mathcal L$:
\begin{gather}\label{t1}
\delta(x+l)=e^{2\pi i\langle\varrho,l\rangle}\delta(x)\,,\\
\label{t22}\delta(x+l\tau)=e^{2\pi i\langle \varrho,l\rangle-2\pi
i\langle \Omega l,x\rangle-\pi i\langle \Omega
l,l\rangle\tau}\delta(x)\,.
\end{gather}
As a corollary of \eqref{char} and \eqref{t1}-\eqref{t22}, we
conclude that the numerator $\Phi$ in \eqref{gans} must have the
translation properties as follows:
\begin{gather*}
\Phi(x+l)=e^{2\pi i\langle a+\varrho,l\rangle}\Phi(x)\,,\\
\Phi(x+l\tau)=e^{2\pi i\langle b+\varrho,l\rangle-2\pi i\langle
\Omega l,x\rangle-\pi i\langle \Omega l,l\rangle\tau}\Phi(x)\,.
\end{gather*}
The vector space of entire functions with such properties is
finite-dimensional. Indeed, let $\omega(x,y)$ denote the bilinear
form on $V$ associated with the operator $\Omega$:
\begin{equation}\label{q}
\omega(x,y)=\sum_{\alpha\in\mathcal
A}m_\alpha\langle\alpha,x\rangle\langle\alpha,y\rangle\,.
\end{equation}
It is symmetric positive and integer-valued on the lattice
$\mathcal L$. Let us define $\Theta\genfrac{[}{]}{0pt}{}{p}{q}$ by
the following series:
\begin{equation}\label{bigtheta}
\Theta\genfrac{[}{]}{0pt}{}{p}{q}(x)=\sum_{l\in\mathcal
L}\exp(2\pi i\omega(l+p,x+q)+\pi i\tau \omega(l+p,l+p))\,.
\end{equation}
It has the following translation properties:
\begin{gather}\label{thtr1}
\Theta\genfrac{[}{]}{0pt}{}{p}{q}(x+l)=e^{2\pi i \omega(l,p)
}\Theta\genfrac{[}{]}{0pt}{}{p}{q}(x)\,,\\\label{thtr2}
\Theta\genfrac{[}{]}{0pt}{}{p}{q}(x+l\tau)=e^{-2\pi i\omega
(l,x+q)-\pi
i\tau\omega(l,l)}\Theta\genfrac{[}{]}{0pt}{}{p}{q}(x)\,.
\end{gather}
It is easy to show (see e.g.\cite{M}) that the space of
holomorphic functions with such translation properties has
dimension equal to $[\mathcal M:\Omega\mathcal L]$; this is equal
to $\det(\Omega_{ij})$ where $\Omega_{ij}=\langle \Omega
e_i,e_j\rangle$ for some basis $e_1,\dots,e_n$ of $\mathcal L$. A
natural basis in this space is given by the functions
$\Theta\genfrac{[}{]}{0pt}{}{p}{q+r}$ with
$r\in\Omega^{-1}\mathcal M$ running over the set of
representatives in $\Omega^{-1}(\mathcal M)/\mathcal L$.

For the later purposes, let us use slightly different basis,
namely, the functions
\begin{equation}\label{basis}
\Phi_r=e^{\langle k,x\rangle}\Theta(x+\gamma+r)\,,\qquad
\Theta:=\Theta\genfrac{[}{]}{0pt}{}{0}{0}\,,\quad
r\in\Omega^{-1}(\mathcal M)/\mathcal L\,.
\end{equation}
It is easy to relate the parameters $k\in V^*$ and $\gamma\in V$
to $p,q$: $$k=2\pi i\Omega p\,,\qquad \gamma=q+\tau p\,.$$ Let us
denote the linear space generated by the functions \eqref{basis}
as $\mathcal U_{k,\gamma}$.

We conclude that $\psi$ must belong to linear space
$\delta^{-1}\mathcal U_{k,\gamma}$ with $k,\gamma$ related in a
simple way to the 'quasimomenta' $a,b$ in \eqref{char}. Now recall
Proposition \ref{f}. It implies that $\psi$ must belong to the
(finite-dimen\-sio\-nal) subspace $\mathcal
W_{k,\gamma}=\left(\delta^{-1}\mathcal U_{k,\gamma}\right)\cap
\mathcal W$. It also implies that $L$ (as well as any of $L_i$'s)
preserves this finite-dimensional space, so we can eventually find
$\psi$ by diagonalizing the action of $L$ on $\mathcal
W_{k,\gamma}$. Note that since the double-Bloch solutions $\psi$
form an $n$-parametric family (with
$\lambda=(\lambda_1,\dots,\lambda_n)$ being the parameters), this
space $\mathcal W_{k,\gamma}$ will be nonzero only for
$(k,\gamma)$ belonging to a certain $n$-dimensional subvariety. In
most cases $\dim\mathcal W_{k,\gamma}\le 1$, so the (unique)
function $\psi_{k,\gamma}$ generating $W_{k,\gamma}\ne 0$ will be
an eigenfunction for $L$ automatically.

All this simplifies a little for the Calogero--Moser models, so
let us consider this case in more detail. For a given reduced,
irreducible root system $R$ in $V^*$ and a fixed $W$-invariant
$m(\alpha)$, we consider the Calogero--Moser operator
\begin{equation}\label{ca}
L=-\Delta+\sum_{\alpha\in
R_+}m_\alpha(m_\alpha+1)(\alpha,\alpha)\wp(\langle\alpha,x\rangle|\tau)\,.
\end{equation}
The Bloch solutions must be of the form $\psi=\delta^{-1}\Phi$ as
in \eqref{gans}. Note that $\delta$ in this case has the following
symmetry:
\begin{equation}\label{wt}
\delta(wx)=\varepsilon_m(w)\delta(x)\quad\text{for any}\ w\in W\,,
\end{equation}
where $\varepsilon_m$ is the one-dimensional character of $W$ such
that
\begin{equation}\label{eps}
\varepsilon_m(s_\alpha)=(-1)^{m_\alpha}\,.
\end{equation}
The bilinear symmetric form \eqref{q} is obviously $W$-invariant,
and since $R$ is irreducible, $\omega$ must be proportional to the
$W$-invariant scalar product $(x,y)$. So,
$\omega(x,y)=\kappa\cdot(x,y)$ for some $\kappa$ which depends on
$m_\alpha$. (For instance, if $R$ consists of one $W$-orbit only
and $m_\alpha\equiv m$, one has $\omega(x,y)=mh(x|y)$ where
$h=h(R)$ is the Coxeter number and the form $(x|y)$ on $V$ is
normalized in such a way that $(\alpha^\vee|\alpha^\vee)=2$ for
all $\alpha^\vee\in R^\vee$.)

The lattice $\mathcal L$ in this case is the coweight lattice
$P^\vee$ of $R$, while $\mathcal M$ is the root lattice $Q$.
Still, the numerator $\Phi$ must belong to the finite-dimensional
space $\mathcal U_{k,\gamma}$ spanned by the functions
\eqref{basis}.

We already know that a Bloch solution $\psi=\delta^{-1}\Phi$
appears only for those $k,\gamma$ when $\delta^{-1}\mathcal
U_{k,\gamma}\cap\mathcal W\ne 0$. Such $k,\gamma$ can be
effectively determined. Indeed, due to \eqref{wt}, the conditions
\eqref{qq} for $\psi$ reduce to the quasi-invariance of $\Phi$:
\begin{equation}
\label{cq} \Phi(x)-\Phi(s_\alpha x)\quad\text{is divisible
by}\quad \langle\alpha,x\rangle^{2m_\alpha+1}\quad\text{for all}\
\alpha\in R\,.
\end{equation}
(These are local conditions near $\pi_\alpha=\{x:
\langle\alpha,x\rangle=0\}$, similar conditions for other
hyperplanes $\pi_\alpha^{m,n}\in Sing$ will follow because $\psi$
is quasiperiodic.)

These conditions can be rewritten as
\begin{equation}\label{der}
\langle{\alpha^\vee},{\partial}\rangle^{2j-1}\Phi\equiv 0
\qquad\text{for}\ \langle\alpha,x\rangle=0\quad\text{and}\
j=1,\dots,m_\alpha\,,
\end{equation}
with $\langle\alpha^\vee,\partial\rangle$ denoting the derivative
in $\alpha^\vee$-direction.

Now recall that we have the period lattice $\mathcal L=P^\vee$,
and $\Phi$ belongs to the linear space $\mathcal U_{k,\gamma}$ of
the functions with the translation properties
\eqref{thtr1}--\eqref{thtr2}. Consider the following sublattice
$\mathcal L^\alpha\subseteq\mathcal L$:
\begin{equation}
\label{ort}
\mathcal L^\alpha:=\mathcal L'\oplus\mathcal
L''\,,\qquad \mathcal L'=\mathcal L\cap\ker\alpha\,,\quad\mathcal
L''=\mathcal L\cap\mathbb R\alpha^\vee\,.
\end{equation}
Let $\mathcal U_{k,\gamma}^\alpha$ denote the space of theta
functions with the same translation properties, but for the
translations $l$ from $\mathcal L^\alpha$ only. Obviously, we have
a natural inclusion map $\mathcal
U_{k,\gamma}\hookrightarrow\mathcal U_{k,\gamma}^\alpha$. It is
possible to describe this linear map explicitly, using the
standard bases in both spaces (look at the formula \eqref{add}
below which is a particular example of such relation). According
to \eqref{ort}, the lattice $\mathcal L^\alpha$ is the direct
orthogonal sum of two sublattices. Thus, the corresponding theta
functions from $\mathcal U_{k,\gamma}^\alpha$ will be the products
of the $(n-1)$-dimensional theta functions related to $\mathcal
L'$ and the one-dimensional theta functions related to the lattice
$\mathcal L''$. This corresponds to the decomposition of $\mathcal
U_{k,\gamma}^\alpha$ into a tensor product: $\mathcal
U_{k,\gamma}^\alpha=(\mathcal U_{k,\gamma}^\alpha)'\otimes
(\mathcal U_{k,\gamma}^\alpha)''$. Applying derivative in
$\alpha^\vee$ direction will affect the one-dimensional theta
functions only. As a result, for each $j$ we have an explicit
linear map $\Gamma^{\alpha,j}$ from $\mathcal U_{k,\gamma}$ to
$(\mathcal U_{k,\gamma}^\alpha)'$ given by
$$\Gamma^{\alpha,j}\varphi=\langle{\alpha^\vee},{\partial}\rangle^
{2j-1}\varphi|_ {\langle\alpha,x\rangle=0}\,.$$ This map is given
by a matrix whose entries are certain combinations of
one-dimensional theta-functions and their derivatives. Now we can
organize all these maps for $\alpha\in R_+$ into one big linear
map
\begin{equation}
\Gamma: \mathcal U_{k,\gamma}\mapsto\bigoplus_
{\genfrac{}{}{0pt}{}{\alpha \in R_+}{j=1 \dots m_\alpha}}
(\mathcal U^\alpha_{k,\gamma})'\,,\qquad
\Gamma=(\Gamma^{\alpha,j})_{\alpha\in R_+, j=1,\dots, m_\alpha}\,,
\end{equation}
with $\Gamma^{\alpha,j}$ defined above. We can think of $\Gamma$
as an $M\times N$ matrix, where $N, M$ are the dimensions of the
source and the target spaces, respectively.

The outcome is the following: a double-Bloch solution $\psi$
appears exactly for those $k,\gamma$ where this linear map has
nontrivial kernel. This gives equations on such $k,\gamma$ (by
equating to zero all $N\times N$ minors of the matrix $\Gamma$).
In its turn, the kernel will determine a corresponding Bloch
eigenfunction. (If the kernel has dimension $>1$, it still defines
an invariant subspace for the action of $L$, so we have at least
one double-Bloch solution.)

So, let $\widetilde {\mathcal C}$ denote an analytic subvariety in
$V^*\times V$ given by
\begin{equation}\label{C}
\widetilde {\mathcal C}=\{(k,\gamma) |\  \ker \Gamma\ne 0\}\,.
\end{equation}
Formulas \eqref{thtr1}-\eqref{basis} make clear that $\widetilde
{\mathcal C}$ is invariant under the following transformations of
$k,\gamma$:
\begin{gather}\label{fac1}
(k,\gamma) \mapsto (k,\gamma+l)\,,\quad l\in \Omega^{-1}\mathcal
M\,,\\\label{fac2} (k, \gamma)\mapsto (k+2\pi i\Omega
l,\gamma+\tau l)\,,\quad l\in \Omega^{-1}\mathcal M\,.
\end{gather}
Now, for a function of one variable $f(z)=e^{kz}\theta(z+c)$ its
derivatives at $z=0$ are obviously polynomial in $k$. Therefore
$\widetilde {\mathcal C}$, after being factored by the
translations above, can be considered as an algebraic covering of
an abelian variety (a product of elliptic curves) $\mathbb
C^n/\mathcal L'+\tau\mathcal L'$ where $\mathcal
L'=\Omega^{-1}\mathcal M$.

\begin{prop}
The double-Bloch eigenfunctions of the Calogero--Moser operator
\eqref{ca} are parametrized by the points of an algebraic variety
which is a covering of an abelian variety $\mathbb C^n/\mathcal
L'+\tau\mathcal L'$ where $\mathcal L'=\Omega^{-1}Q$, $Q$ is the
root lattice and the map $\Omega$ is defined by \eqref{Om}.
\end{prop}

A similar analysis applies to any (integrable) generalized Lam\'e
operator, so the double-Bloch eigenfunctions are also parametrized
by the points of an algebraic variety covering a product of
elliptic curves.

\medskip

Now let $\mathcal C$ denote the result of factoring the variety
\eqref{C} by the translations \eqref{fac1}--\eqref{fac2}. It is an
algebraic variety parametrizing the double-Bloch eigenfunctions of
$L$. Below, following \cite{FV}, we will refer to it as the {\bf
Hermite--Bloch variety} for $L$. It differs from the complex Bloch
variety, traditionally defined as the set of $(\mu,E)\in (\mathbb
C^\times)^n\times\mathbb C$ such that there exists $\psi$ with
$L\psi=E\psi$ and $\psi(x+l_i)=\mu_i\psi(x)$ where $l_1,\dots l_n$
is a basis in $\mathcal L$. Note that the latter is a
transcendental complex analytic variety.

\medskip
\begin{remark}\label{winv}
In case of the Calogero--Moser operator related to a root system
$R$, there is a natural action of the Weyl group $W$ on the
Hermite--Bloch variety $\mathcal C$. Also, there is a natural
projection of $\mathcal C$ onto $\mathbb C^n$ sending
$\psi=\psi_{k,\gamma}$ to the set of eigenvalues $\lambda_i$,
$L_i\psi=\lambda_i\psi$. This is a $|W|$-sheeted covering and the
Weyl group acts on $\mathcal C$ by permuting the points in the
fiber.
\end{remark}

\begin{remark} Note that our results do not contradict the theorem of
Feldman--Kn\"orrer--Trubowitz \cite{FKT} in dimension two, since
their result only applies to a {\it real-valued smooth} potential
$u$ in $\mathbb R^2$.
\end{remark}

\medskip
The Hermite--Bloch variety $\mathcal C$ is a subvariety in the
total space of a certain bundle over the product of elliptic
curves defined by \eqref{fac1}--\eqref{fac2}. This bundle
naturally compactifies to a bundle with the fibers isomorphic to
the projective space $\mathbb P^n$. As a result, $\mathcal C$
compactifies to a projective variety, covering the product of
elliptic curves.

The variety $\mathcal C$ is closely related to the so-called {\bf
spectral variety} which is defined as follows. Suppose $L$ is a
strongly integrable generalized Lam\'e operator, so we have $n$
commuting operators $L_1=L,\dots, L_n$, which generate a
commutative subalgebra in the ring of PDO with meromorphic
elliptic coefficients. Then by \cite{BEG}, theorem 2.2, the
centralizer of this subalgebra will be a maximal commutative ring
which we will denote by $\mathcal Z(L)$ (using Proposition \ref{f}
one can show that the operators in this ring will share a common
family of the double-Bloch eigenfunctions). Each operator in
$\mathcal Z(L)$ must have constant highest symbols, by Lemma
\ref{ber}. Then from the strong integrability we immediately
derive that $\mathcal Z(L)$ is finitely generated. Thus,
$\mathrm{Spec}\mathcal Z(L)$ defines an affine algebraic variety,
which we call the spectral variety. It is not quite clear whether
the spectral variety is isomorphic to the Hermite--Bloch variety
(for instance, the latter may not be affine), but at least they
must be birationally equivalent.

Finally, let us remark on some algebraic geometry behind the
double-Bloch solutions and the Hermite--Bloch variety for the
Calogero--Moser system \eqref{ca}. We consider the torus
$T=\mathbb C^n/Q^\vee+\tau Q^\vee$ where $Q^\vee=R^\vee\otimes
\mathbb Z$ is the coroot lattice. Let us define the following
subsheaf $\mathcal Q \subset\mathcal O(T)$ of the structure sheaf
of $T$ by requiring its local sections to have zero normal
derivatives of order $1,3,\dots,2m_\alpha-1$ along each of the
hyperplanes $\langle \alpha, x\rangle=0$ (considered as hypertori
in $T$). The sheaf $\mathcal Q$ can be considered as the structure
sheaf $\mathcal O(X)$ of a singular variety $X$, with $T$ being
its injective normalization (cf. \cite{BEG1}). Such $X$ is
projective; it is a $|W|$-sheeted covering of the weighted
projective space $T/W$ considered by Looijenga \cite{Lo}, see also
\cite{BS}. Notice that from the results of \cite{EG} it follows
that $X$ is Cohen-Macaulay and Gorenstein. Let us consider now the
group $\mathrm{Pic}(X)$ of invertible sheaves on $X$. Then each of
the double-Bloch solutions $\psi_{k,\gamma}$ represents a
meromorphic section of a degree zero line bundle on $X$ (to
define degree, we use the pull-back to the torus $T$). In this way
the Hermite--Bloch variety for the Calogero--Moser system becomes
an $n$-dimensional subvariety in $\mathrm{Pic}^0(X)$. It would be
interesting to study this relation in more detail.

\begin{remark} An interesting thing is to analyse how the spectral variety
changes when $\tau$ goes to $+i\infty$ (trigonometric limit). In
this limit the spectral variety becomes rational and is relatively
well understood. Thus, one could think of the whole family
depending on $\tau$ as a deformation of this rational variety.
This point of view was used in \cite{VS} to construct the spectral
surface in the simplest $A_2$ case.
\end{remark}

\subsection{Discrete spectrum eigenstates}
\label{spe}

Let us explain how the Bloch solutions can be used to construct
the discrete spectrum eigenstates of $L$. Our discussion is
strictly confined to the Calogero--Moser operator \eqref{ca}. We
take a purely imaginary $\tau$, this ensures that the potential in
\eqref{ca} is real-valued for $x\in\mathbb R^n$. The
Calogero--Moser operator $L$ is defined on a dense subset of
$L^2(\mathbb R^n)$ and it is self-adjoint only formally, and its
Bloch solutions are singular. It has square-integrable
eigenstates, though. Namely, let $\psi=\psi_{k,\gamma}$ be one of
the double-Bloch solutions constructed in the previous section.
Given such a $\psi$, let us symmetrize it as follows:
\begin{equation}\label{spect}
\Psi(x)=\sum_{w\in W}\varepsilon_m(w)\det w\, \psi(wx)\,,
\end{equation}
where $W$ is the Weyl group of the root system $R$ and
$\varepsilon_m$ is the character \eqref{eps}. The Calogero--Moser
operator $L$ is $W$-invariant, thus the constructed $\Psi$ will be
again its eigenfunction (by the same reason, it will be an
eigenfunction for all commuting operators $L_i$). A priori, $\Psi$
might have poles in $\mathbb R^n$ along the hyperplanes
$\langle\alpha,x\rangle=c$, $c\in\mathbb Z$. However, it is easy
to see that $\Psi$ has no poles along the hyperplanes
$\langle\alpha,x\rangle=0$. This follows immediately from the
properties \eqref{qq} of $\psi$. To avoid the appearance of
singularities on other hyperplanes $\langle\alpha,x\rangle=c$, one
has to impose the condition that all the terms $\psi(wx)$ in the
sum \eqref{spect} have the same Bloch--Floquet multipliers with
respect to a shift $x\mapsto x+l$ with $l\in\mathcal L=P^\vee$.
This means that $\exp\langle k,l\rangle=\exp\langle wk,l\rangle$
for all $w\in W$ and $l\in P^\vee$, which in its turn implies that
$k$ belongs to the lattice $2\pi i P$. So, we have the following
result.

\begin{prop}\label{sp1} Let $L$ be the Calogero--Moser operator \eqref{ca}.
Then for any point $(k,\gamma)$ of its Bloch--Hermite variety
which satisfies an additional condition $k\in 2\pi i P$ (with $P$
being the weight lattice for $R$), the corresponding function
\eqref{spect} (if nonzero) will be a nonsingular in $\mathbb R^n$
eigenfunction of the Calogero--Moser operator \eqref{ca} and of
the higher operators $L_2,\dots, L_n$.
\end{prop}

By construction, $\Psi$ vanishes along the hyperplanes
$\langle\alpha,x\rangle\in\mathbb Z$, and it gets a factor of
$(-1)^{m_\alpha+1}$ under the orthogonal reflection with respect
to such a hyperplane. Since these are the reflection hyperplanes
of the affine Weyl group of $R$, they cut $\mathbb R^n$ into its
fundamental domains (alcoves), so the restriction of $\Psi$ to
each alcove will be, essentially, the same. We can restrict $\Psi$
to any alcove, extending it by zero outside, and this gives us a
finitely supported smooth eigenfunction of $L$ (notice that in the
complex domain it still has poles). We see from this that the
discrete spectrum of $L$ in $L^2(\mathbb R^n)$ is infinitely
degenerate (one says that the spectral problem for $L$ splits into
identical spectral problems on each of the alcoves). Morally, this
is the reason why one should expect the same spectrum considering
$L$ not on $L^2(\mathbb R^n)$ but on the space $L^2(T)^W$ of
$W$-invariant functions on the torus $T=\mathbb R^n/Q^\vee$, as in
\cite{KT}. The latter case is simpler from the technical point of
view, since the operator $L$ is essentially self-adjoint on
$L^2(T)^W$, see \cite{KT} for the details.

In \cite{KT} Komori and Takemura considered the elliptic
Calogero--Moser problems as a perturbation (in $\tau$) of the
trigonometric case $\tau=+i\infty$, and Theorem 3.7 of \cite{KT}
claims that for sufficiently small $p=e^{2\pi i\tau}$ the family
of eigenfunctions (Jack polynomials) which corresponds to $p=0$,
admits analytic continuation in $p$, and the resulting functions
will give rise to a complete orthogonal family of eigenfunctions
of $L$ in $L^2(T)^W$. One can show that our family in the limit
$\tau\to +i\infty$ specializes to the Jack polynomials. Comparing
this with the previous discussion, we conclude that our family
must coincide with the one considered in \cite{KT}.

\section{Calogero--Moser model of $B_2$ type}

In this section we consider the following 2-dimensional
Schr\"odinger operator
\begin{equation}\label{b2}
L=-\Delta+2\wp(x_1)+2\wp(x_2)+4\wp(x_1-x_2)+4\wp(x_1+x_2)\,,
\end{equation}
where $\wp(z)=\wp(z|\tau)$ is the Weierstrass $\wp$-function with
the periods $1,\,\tau$ ($\mathrm{Im}\,\tau > 0$). Our goal is to
calculate its double-Bloch eigenfunctions, i.e. such $\psi$ that
\begin{align}\label{rebl}
\psi(x+e_j)&=\lambda_j\psi(x)\\ \label{imbl}\psi(x+\tau
e_j)&=\mu_j\psi(x)\quad (j=1,2)\,,
\end{align}
where $(e_1, e_2)$ is the standard basis in $\mathbb C^2$ and
$(\lambda_1, \lambda_2, \mu_1, \mu_2)$ are fixed Bloch--Floquet
multipliers.

First we recall some standard definitions and formulas from the
theory of theta-functions, see \cite{BE,M}. Let
$\theta\genfrac{[}{]}{0pt}{}{\alpha}{\beta}$ be the
one-dimensional theta-function (with characteristics), defined by
the following series: $$
\theta\genfrac{[}{]}{0pt}{}{\alpha}{\beta}(z|\tau)=\sum_{n\in{\mathbb
Z}}\exp(\pi i (n+\alpha)^2\tau+2\pi i(n+\alpha)(z+\beta))\,. $$
Notice that $\alpha$ and $\beta$ are defined modulo 1: $$
\theta\genfrac{[}{]}{0pt}{}{\alpha+1}{\beta}
=\theta\genfrac{[}{]}{0pt}{}{\alpha}{\beta}\,,\quad
\theta\genfrac{[}{]}{0pt}{}{\alpha}{\beta+1}=e^{2\pi
i\alpha}\theta\genfrac{[}{]}{0pt}{}{\alpha}{\beta}\,.$$ Later we
will need the following formula which can be easily derived from
the definitions:
\begin{multline}\label{add}
\theta\genfrac{[}{]}{0pt}{}{\alpha_1}{0}(x_1|\tau)
\theta\genfrac{[}{]}{0pt}{}{\alpha_2}{0}(x_2|\tau)=
\theta\genfrac{[}{]}{0pt}{}{\alpha_+}{0}(x_+|2\tau)
\theta\genfrac{[}{]}{0pt}{}{\alpha_-}{0}(x_-|2\tau)\\+
\theta\genfrac{[}{]}{0pt}{}{\alpha_++\frac12}{0}(x_+|2\tau)
\theta\genfrac{[}{]}{0pt}{}{\alpha_-+\frac12}{0}(x_-|2\tau)\,,
\end{multline}
where $\alpha_{\pm}=\frac12(\alpha_1\pm\alpha_2),\ x_\pm=x_1\pm
x_2$.

We will mostly use $\theta\genfrac{[}{]}{0pt}{}{1/2}{1/2}$ which
we will denote simply by $\theta(z)$, which will always stand for
the odd Jacobi theta function \eqref{theta1}.

According to the previous section, $\psi$ must have the form
\begin{equation}\label{ansatz}
\psi = \frac{\Phi(x_1,x_2)}
{\theta(x_1|\tau)\theta(x_2|\tau)\theta(x_1-x_2|\tau)\theta(x_1+x_2|\tau)}\,.
\end{equation}
Here $\Phi$ is nonsingular in $\mathbb C^2$. The translation
properties for $\psi$ easily translate into the properties of
$\Phi$:
\begin{gather*}
\Phi(x+e_j)=-\lambda_j\Phi(x)\,,\\ \Phi(x+\tau e_j)=-\mu_je^{-3\pi
i\tau-6\pi ix_j}\Phi(x)\,.
\end{gather*}
Standard considerations from the theory of theta-functions show
that the linear space of functions with these properties has
dimension $9$ and $\Phi$ must be of the form
\begin{equation}\label{phi}
\Phi=\exp(K_1x_1+K_2x_2)\sum_{0\le i,j\le
2}c_{ij}\theta\genfrac{[}{]}{0pt}{}{i/3}{0}(3x_1+\gamma_1|3\tau)\theta\genfrac{[}{]}{0pt}{}{
j/3}{0}(3x_2+\gamma_2|3\tau)\,,
\end{equation}
where $c_{ij}$ are arbitrary constants and parameters $\gamma_j,
K_j$ relate to $\lambda_j, \mu_j$ as follows:
\begin{equation}\label{1.4}
\lambda_j=-e^{K_j}\,,\quad \mu_j=-e^{-2\pi i\gamma_j+K_j\tau}\,.
\end{equation}

\bigskip
\begin{remark}\label{rem1}
The shifting
\begin{equation}\label{a}
\gamma_j\mapsto \gamma_j+1
\end{equation}
does not change the space \eqref{phi}, the same applies to the
shifts
\begin{equation}\label{ak}
(K_j, \gamma_j)\mapsto (K_j+2\pi i, \gamma_j+\tau)\,.
\end{equation}
Conversely, for any given $\lambda_j,\mu_j$ the corresponding
$(\gamma_j,K_j)$ are determined uniquely modulo shifts
\eqref{a}--\eqref{ak}.
\end{remark}

\bigskip
Now, in accordance with Proposition \ref{f}, we impose the
following 'vanishing' conditions on $\Phi$:
\begin{align}\label{ax1}
  \d_1\Phi \equiv 0 & \quad \text{for}\ x_1=0\,, \\ \label{ax2}
  \d_2\Phi \equiv 0 & \quad \text{for}\ x_2=0\,, \\ \label{ax+}
  (\d_1+\d_2)\Phi \equiv 0 & \quad \text{for}\ x_1+x_2=0\,, \\
  \label{ax-}
  (\d_1-\d_2)\Phi \equiv 0 & \quad \text{for}\ x_1-x_2=0\,.
\end{align}
As we will see below, for a certain $2$-dimensional surface in
$4$-dimensional space of parameters $(k_j, a_j)$, the conditions
\eqref{ax1}--\eqref{ax-} cut a one-dimensional subspace in
$9$-dimensional space \eqref{phi}. Thus, the corresponding $\psi$
will be an eigenfunction for $L$ automatically.

To determine the corresponding $(K_j, \gamma_j)$, let us rewrite
$\Phi$ using \eqref{add} and making identification
$\theta\genfrac{[}{]}{0pt}{}{\alpha
+1}{0}=\theta\genfrac{[}{]}{0pt}{}{\alpha}{0}$:
\begin{multline*}
\Phi= e^{K_+x_+ + K_-x_-}\sum_{0\le l,m \le 2} \widetilde
c_{lm}\left\{\theta\genfrac{[}{]}{0pt}{}{ l/3}{0}(3x_+ +
\gamma_+|6\tau)\theta\genfrac{[}{]}{0pt}{}{ m/3}{0}(3x_- +
\gamma_-|6\tau)\right.
\\ \left. + \theta\genfrac{[}{]}{0pt}{}{ l/3+1/2}{0}(3x_+ +
\gamma_+|6\tau)\theta\genfrac{[}{]}{0pt}{}{ m/3+1/2}{0}(3x_- +
\gamma_-|6\tau)\right\} \,,
\end{multline*}
where $K_\pm=\frac 12(K_1\pm K_2)$,\ $\gamma_\pm=\gamma_1\pm
\gamma_2$ and
\begin{equation}\label{tildec}
\widetilde c_{lm}=c_{ij}\quad \text{with}\quad i\equiv
l+m\pmod{3}\,,\ j\equiv l-m\pmod{3}\,.
\end{equation}
It is easy to see now that \eqref{ax+} leads to six linear
equations on $\widetilde c_{lm}$ of the form
\begin{equation*}
\sum_{l=0}^2 A_{l}\widetilde c_{lm}=0\,,\qquad \sum_{l=0}^2
B_{l}\widetilde c_{lm}=0\qquad(m=0,1,2)
\end{equation*}
for certain explicitly given $A=(A_l)$, $B=(B_l)$. For generic
parameters $\gamma_1,\gamma_2$ the vectors $A,B$ will be linearly
independent. Therefore, these $6$ equations determine the
$2$-dimensional kernel of the matrix $\widetilde C = (\widetilde
c_{lm})$. Similarly, \eqref{ax-} gives six more equations for
$\widetilde c_{lm}$, which determine the cokernel of $\widetilde
C$. Thus, for any $K_1,K_2$ and generic $\gamma_1,\gamma_2$ the
vanishing conditions \eqref{ax+}--\eqref{ax-} determine
$\widetilde c_{lm}$ (and, hence, $\Phi$) uniquely up to a common
factor. In principle, it is straightforward to write down explicit
expressions for the coefficients $c_{ij}$ but they are cumbersome
and not very useful. However, there is a better way of getting an
expression for $\Phi$, by taking a limit $\o\to 0$ in the formula
for the difference case, see \eqref{qans} below. It turns out that
$\Phi$ has the form
\begin{equation}\label{up}
\Phi=\frac{e^{(k,x)}}{\t(a_1|\tau)\t(a_2|\tau)} \sum_{0\le i,j\le
2} b_{ij}(k_1+k_2)^i(k_1-k_2)^j \,,
\end{equation}
where the parameters $k_1,k_2$ and $a_1,a_2$ relate to
$K_j,\gamma_j$ as
\begin{equation}\label{K}
k_j=K_j+\pi i\,,\qquad a_j=\gamma_j-(1+\tau)/2\,.
\end{equation}
The coefficients $b_{ij}$ depend on $x$ and $a_1,a_2$ and are
given by the following recipe. Let us introduce formal commutative
variables $A,B,C,D$. We also need the following scalar $\lambda$
depending on $\tau$:
$$\lambda:=\theta'''(0|\tau)/\theta'(0|\tau)\,.$$ Now introduce
$U,V$ as $U:=A+B-C$,\quad$V:=A+B-D$ and put
\begin{align}\notag b_{22}=&1\,,\quad b_{21}=2U\,,\quad
b_{12}=2V\,,\\\label{upp} b_{20}=&-\lambda+U^2\,,\quad
b_{02}=-\lambda+V^2\,,\quad b_{11}=4UV\,,
\\\notag b_{10}=&2V(-\lambda+U^2)\,,\quad
b_{01}=2U(-\lambda+V^2)\,,\\\notag
b_{00}=&(-\lambda+U^2)(-\lambda+V^2)\,.
\end{align}
After that one opens the brackets, so each of $b_{ij}$ becomes a
sum of monomials in $A,B,C,D$ with scalar coefficients, and then
replaces each monomial using the following rule:
\begin{equation}\label{uppe}
cA^pB^qC^rD^s\longrightarrow
c\t^{(p)}(x_1+a_1)\t^{(q)}(x_2+a_2)\t^{(r)}(x_1-x_2)
\t^{(s)}(x_1+x_2)\,.
\end{equation}
Here $\t=\t(z|\tau)$, as before, denotes the odd Jacobi theta
function, and the upper index in brackets refers to taking
derivatives in $z$. We treat a scalar as a multiple of
$A^0B^0C^0D^0$ assuming, as usual, that $f^{(0)}=f$. To illustrate
this, we present below some first of the coefficients:
\begin{align}\notag
b_{22}=&\t(x_1+a_1)\t(x_2+a_2)\t(x_1-x_2)\t(x_1+x_2)\,,\\ \notag
b_{21}=&2\t'(x_1+a_1)\t(x_2+a_2)\t(x_1-x_2)\t(x_1+x_2)-\\
\label{upper}
&2\t(x_1+a_1)\t'(x_2+a_2){\theta}(x_1-x_2){\theta}(x_1+x_2)-\\
\notag&
2\t(x_1+a_1)\t(x_2+a_2)\t'_1(x_1-x_2){\theta}(x_1+x_2)\,,\\\notag
b_{20}=&
-\lambda\t(x_1+a_1)\t(x_2+a_2){\theta}(x_1-x_2){\theta}(x_1+x_2)+\dots
\end{align}

\begin{prop} For any $K_1,K_2$ and generic $\gamma_1,\gamma_2$
there is a unique (up to a constant factor) function $\Phi(x)$ of
the form \eqref{phi} with the properties \eqref{ax1}--\eqref{ax-}.
It is described by the formulas \eqref{up}--\eqref{uppe} above.
\end{prop}

To prove the formula \eqref{up}, one goes to the limit $\o\to 0$
in the formula \eqref{qans} from the next section, picking up the
first nonzero term (of order $6$ in $\o$).

 \medskip
Formulas \eqref{up}--\eqref{uppe} fix the dependence of $\Phi(x)$
on $4$ parameters $a_1,a_2,k_1,k_2$. Let us consider now the
translation properties of $\Phi$ regarded as a function of these
parameters. Recall that for generic $a_1, a_2$ the function $\Phi$
was determined uniquely up to a factor by \eqref{phi} and
\eqref{ax+}--\eqref{ax-}. Thus, it follows immediately from Remark
\ref{rem1} that under the shifts \eqref{a}--\eqref{ak} $\Phi$ must
remain the same, up to a factor independent on $x$. To find this
factor, it is sufficient to look at the formula \eqref{upper} for
the leading coefficient $b_{22}$. As a result, we conclude that
the function $\Phi$, given by formulas \eqref{phi},
\eqref{up}--\eqref{uppe} is invariant with respect to the shifts
\eqref{a}--\eqref{ak}:
\begin{gather}\label{t12}
\Phi(a_1+1, a_2, k_1, k_2)=\Phi(a_1, a_2+1, k_1, k_2)=\Phi(a_1,
a_2, k_1, k_2)\,,\\\label{t2} \Phi(a_1+\tau, a_2, k_1+2\pi i,
k_2)=\Phi(a_1, a_2+\tau, k_1, k_2+2\pi i)=\Phi(a_1, a_2, k_1,
k_2)\,.
\end{gather}

Now let us find the interrelations between the parameters $a_j,
k_j$ which will guarantee two remaining vanishing conditions
\eqref{ax1}--\eqref{ax2}. Let $G_1$ and $G_2$ denote the
derivatives $\d_1\d_2\Phi$ and $\d_1\d_2^3\Phi$ evaluated at
$x_1=x_2=0$:
\begin{equation}\label{G}
G_1=\frac{\d^2 \Phi}{\d x_1\d x_2}(0,0)\,,\qquad G_2=\frac{\d^4
\Phi}{\d x _1\d x_2^3}(0,0)\,.
\end{equation}
Thus defined $G_1, G_2$ will be regarded as functions of the
parameters $a_1,a_2,k_1,k_2$. Consider now the following two
equations on these $4$ parameters:
\begin{equation}\label{ob}
G_1=0\,,\quad G_2=0\,.
\end{equation}
It is clear that conditions \eqref{ax1} imply both of the
equations \eqref{ob}. Indeed, they guarantee that the function
\begin{equation}\label{fu}
f(t)=\d_1\Phi(0,t)
\end{equation}
is identically zero, in particular, $f'(0)=G_1=0$ and
$f'''(0)=G_2=0$.

More interestingly, \eqref{ob} are 'almost' equivalent to
\eqref{ax1}. To see this, note that the function \eqref{fu} has
the following translation properties in $t$:
\begin{gather*}
f(t+1)=e^{k_2}f(t)\,,\\f(t+\tau)=e^{k_2\tau-2\pi i(3t+a_2)-3\pi
i\tau}f(t)\,.
\end{gather*}
{From} \eqref{ax+}--\eqref{ax-} we know that
$$(\d_1\pm\d_2)\Phi(0,0)=0\,.$$ This gives that
$f(0)=\d_1\Phi(0,0)=0$, while the first equation in \eqref{ob}
gives that $f'(0)=0$. Together with the translation properties
above this implies that $f$ is proportional to the following theta
function: $$e^{k_2t}\theta(t+a_2|\tau)[\theta(t|\tau)]^2\,.$$ This
function has nonzero third derivative at $t=0$ as soon as
\begin{equation}\label{gen}
k_2\theta(a_2|\tau)+\theta'(a_2|\tau)\neq 0\,.
\end{equation}
Thus, the conditions $f(0)=f'(0)=f'''(0)=0$ imply that $f$ is
identically zero provided that \eqref{gen} is true. The outcome is
the following: for all $a_1, a_2, k_1, k_2$ satisfying the
condition \eqref{gen} the equations \eqref{ob} imply \eqref{ax1}.

Our next remark is that for any $a_1, a_2, k_1, k_2$ we have the
relation $\d_1\d_2^3\Phi(0,0)=\d_1^3\d_2\Phi(0,0)$. This is due to
the identity $$ 4\d_1^3\d_2\Phi-4\d_1\d_2^3\Phi=
(\d_1+\d_2)^3(\d_1-\d_2)\Phi- (\d_1-\d_2)^3(\d_1+\d_2)\Phi$$ where
the right-hand side vanishes at $x_1=x_2=0$ due to
\eqref{ax+}--\eqref{ax-}. Thus, repeating the same arguments as
above, we conclude that the equations \eqref{ob} imply also
\eqref{ax2}, as soon as
\begin{equation}\label{gen2}
k_1\theta(a_1|\tau)+\theta'(a_1|\tau)\neq 0\,.
\end{equation}

As a result, we see that under assumptions
\eqref{gen}--\eqref{gen2} both vanishing conditions
\eqref{ax1}--\eqref{ax2} are equivalent to one system \eqref{ob}.
To remove the restrictions \eqref{gen}--\eqref{gen2}, let us
consider equations \eqref{ob} in more details.

First notice that the functions $G_1,G_2$ share the same
translation properties \eqref{t12}--\eqref{t2} in
$(a_1,a_2,k_1,k_2)$ with $\Phi$ (since differentiating in $x$
doesn't affect these properties). Notice also that $G_1, G_2$ are
polynomials in $k_1,k_2$. Let us pass from $k_1,k_2$ to another
variables $p_1,p_2$ as follows:
\begin{equation}\label{p}
p_1=k_1+\zeta(a_1)\,,\quad p_2=k_2+\zeta(a_2)\,,
\end{equation}
where we used $\zeta=\zeta(z|\tau)$ to denote the logarithmic
derivative of $\theta(z)$:
$$\zeta(z)=\frac{\theta'(z|\tau)}{\theta(z|\tau)}\,.$$ This is
slightly different from the Weierstrass $\zeta$-function.

Clearly, $G_1,G_2$ are still polynomials in $p_1,p_2$ with the
coefficients depending on $a_1,a_2$. The translation properties
\eqref{t12}--\eqref{t2} imply that the coefficients in these
polynomials will be elliptic functions of $a_1$ and $a_2$. In
fact, one can write down $G_1, G_2$ quite explicitly. The
following result follows from our calculations for the difference
case from the next section.

\begin{prop}\label{b2var} The system $G_1=G_2=0$ is equivalent to the following system:
\begin{equation}\label{var}
\begin{array}l
p_1(p_2^3+3\zeta'_2p_2+\zeta''_2)=p_2(p_1^3+3\zeta'_1p_1
+\zeta''_1)\,,\\ p_1(p_2^5+10\zeta'_2p_2^3+10\zeta''_2p_2^2+
(5\zeta'''_2+15(\zeta'_2)^2)p_2+\zeta''''_2+10\zeta'_2\zeta''_2)\\=
p_2(p_1^5+10\zeta'_1p_1^3+10\zeta''_1p_1^2+
(5\zeta'''_1+15(\zeta'_1)^2)p_1+\zeta''''_1+10\zeta'_1\zeta''_1)\,,
\end{array}
\end{equation}
where $\zeta_1, \zeta_2$ stand for $\zeta(a_1)$ and $\zeta(a_2)$
while the prime denotes taking the derivative with respect to the
corresponding variable. (Notice that $\wp(z)=-\zeta'+{\rm const}$,
so $\zeta''=-\wp'$ and so on.)
\end{prop}

{From} the discussion below will follow that for generic $a_1,a_2$
the system \eqref{var} has a finite number of solutions
$(p_1,p_2)$. Thus, we can think of \eqref{var} as a finite
covering of the product $\mathcal E\times \mathcal E$ of two
copies of an elliptic curve $\mathcal E=\mathcal E_\tau=\mathbb C
/(\mathbb Z+\tau\mathbb Z)$. In fact, the only $(a_1,a_2)$ where
the fiber is infinite are those with $a_1=\pm a_2\ ({\rm mod}\
1,\tau)$. This corresponds to the following 'vertical' components
of \eqref{var}:
\begin{equation}\label{vert}
\{a_1=a_2\,,\ p_1=p_2\}\,,\qquad\{a_1=-a_2\,,\ p_1=-p_2\}\,.
\end{equation}
%(Note that the function $\Phi$ \eqref{phi} is identically zero
%along these components, and they don't contribute to the Bloch
%variety of $L$.)

Another 'trivial' component is, obviously,
\begin{equation}\label{triv}
\{p_1=p_2=0\}\,.
\end{equation}
If we delete these three components from \eqref{var}, the
remaining part will be, in fact, a $13$-fold covering of $\dot
{\mathcal E}\times \dot {\mathcal E}$. Since we deleted the
component \eqref{triv}, the conditions \eqref{gen}--\eqref{gen2}
and, hence, \eqref{ax1}--\eqref{ax2} are satisfied on the
remaining part. Thus, we arrive at the following theorem.

\begin{theorem}\label{blo} Let $\mathcal C$ be the finite covering of the product
of two (punctured) elliptic curves which is obtained from
\eqref{var} by deleting the components \eqref{vert} and
\eqref{triv}. Then $\mathcal C$ is the Hermite--Bloch variety for
the operator \eqref{b2} and a double-Bloch solution $\psi(x)$
corresponding to a point $(a_1,a_2,p_1,p_2)$ in $\mathcal C$ is
given by the formulas \eqref{ansatz}, \eqref{up}--\eqref{uppe} and
\eqref{p}.
\end{theorem}

We still have to explain why $\mathcal C$ is a $13$-fold covering.
To this end let us consider a family of plane rational curves
$\varphi: \mathbb P^1\to\mathbb P^2$ of degree $5$, depending on a
parameter $a\in \mathcal E=\mathbb C/\mathbb Z+\tau\mathbb Z$ and
defined as follows: if $a\in \mathcal E$ and
$u=(u_0:u_1)\in\mathbb P^1$ then
\begin{align*}
&\varphi(a,u)=(\varphi_0:\varphi_1:\varphi_2)\,,\quad\text{where}\\
&\varphi_0=u_1u_0^4\,,
\\&\varphi_1=u_1^3u_0^2+3\zeta'(a)u_1u_0^4+\zeta''(a)u_0^5\,,
\\&\varphi_2=u_1^5+10\zeta'(a)u_1^3u_0^2+
(5\zeta'''(a)+15(\zeta'(a))^2)u_1u_0^4+
(\zeta''''(a)+10\zeta'(a)\zeta''(a))u_0^5\,.
\end{align*}
Then the solutions $(p_1,p_2)$ of \eqref{var} correspond to the
intersection points of two curves $C_1=\varphi(a_1,\cdot\,)$,
$C_2=\varphi(a_2,\cdot\,)$ from our family. Namely, if
$\varphi(a_1,u)=\varphi(a_2,v)$ then $p_1=u_1/u_0$ and
$p_2=v_1/v_0$ obviously satisfy \eqref{var} and vice versa,
provided $p_1,p_2\ne 0$. We should, however, exclude from
consideration points with $p_1,p_2=\infty$. Namely, all the curves
from our family pass through $(0:0:1)=\varphi(a,\infty)$. It is a
standard exercise in basic algebraic geometry to show that ${\rm
mult}\,(C_1\cap C_2)$ at this point equals $12$. In doing this
local analysis, one immediately observes that the condition
$\zeta(a_1)=\zeta(a_2)$ is necessary for $C_1$ and $C_2$ to
coincide near the point $(0:0:1)$. This proves that the covering
\eqref{var} is finite apart from the components \eqref{vert}.
After that Bezout's theorem tells us that the number of common
points, apart from $(0:0:1)$, equals $5\times 5-12=13$.

\medskip

\begin{remark}\label{weyl} Operator \eqref{b2} is symmetric under the
Weyl group $W$ of $B_2$, whose $8$ elements act by permuting
coordinates and/or changing their signs. This induces the action
of $W$ on Bloch solutions and, therefore, on the Hermite--Bloch
variety \eqref{var}. This action, in terms of $(a_1,a_2,p_1,p_2)$,
is generated by involutions $(a_1,a_2,p_1,p_2)\to
(a_2,a_1,p_2,p_1)$ and $(a_1,a_2,p_1,p_2)\to (-a_1,a_2,-p_1,p_2)$.
\end{remark}

\subsection{Algebraic integrability}

The operator $L$ \eqref{b2} is completely integrable. According to
Theorem \ref{ch}, it has a commuting operator of order four,
$$L_1=\partial_1^2\partial^2_2+\dots\,.$$ Using Proposition
\ref{f}, we see that our $\psi$ is a common eigenfunction for
$L,\, L_1$:
\begin{equation}\label{eig3}
L\psi=E\psi\,,\quad L_1\psi=E_1\psi\,,\quad L_2\psi=E_2\psi\,.
\end{equation}
Here $E, E_1$ are some functions of the parameters $a_1, a_2, k_1,
k_2$ which, in principle, can be calculated explicitly (though we
didn't have enough energy to perform such a calculation).

In \cite{O} it was shown that apart from $L_1$, there is another
operator $$L_2=\partial_1^5-5\partial_1^3\partial_2^2+\dots$$
which commutes with $L, L_1$ (see \cite{O} for the explicit
expression for $L_1,L_2$). The existence of a fifth order quantum
integral $L_2$ means in this case that the Schr\"odinger operator
$L$ is algebraically integrable. To see this directly, let us
consider one more operator $L_3$, obtained from $L_2$ by
interchanging $x_1$ and $x_2$. Then one easily checks that the
common eigenspace \eqref{eig3} of $L$ and $L_1$ is
$8$-dimensional, and for generic $E$ and $E_1$ it is spanned by
the double-Bloch solutions $\psi$, constructed previously. On the
other hand, by proposition \ref{f}, each $\psi$ will be an
eigenfunction of $L_2$ and $L_3$ as well. So, the only thing to
check is that the eigenvalues $E_2, E_3$ separate all $8$
solutions of the system \eqref{eig3}. This is enough to check in
the trigonometric limit $\tau\to+i\infty$, which is not difficult.

\begin{remark} We do not give the precise relation between the Hermite--Bloch
variety and the spectral surface. To find such a relation, a
careful analysis of the structure of the divisor at infinity is
needed. Let us remark that in \cite{Kh} two algebraic relations
between the $4$ operators $L,\dots ,L_3$ were calculated
explicitly. Thus, they determine a $2$-dimensional affine
algebraic variety in $\mathbb C^4$. However, it is not isomorphic
to the spectral surface as an affine variety (though they are
birationally equivalent). This can be seen already in the
trigonometric limit $\tau\to +i\infty$, by using the information
about the spectral variety from \cite{VSC}. Namely, the results of
\cite{VSC} imply that these four operators $L,\dots, L_3$ are not
enough to generate the whole commutative ring (which is isomorphic
to the coordinate ring of the spectral surface).
\end{remark}

\subsection{Spectrum of $L$}
\label{spectrum}

Throughout this subsection we assume that the parameter $\tau$ is
pure imaginary, so the potential $u(x)$ of the Schr\"odinger
operator \eqref{b2} is real-valued for $x\in\mathbb R^2$. Its
singularities is the family of lines
\begin{equation}\label{sing}
x_1\in \mathbb Z\,,\quad x_2\in \mathbb Z\,,\quad x_1+x_2\in
\mathbb Z\,,\quad x_1-x_2\in \mathbb Z\,.
\end{equation}
These lines cut $\mathbb R^2$ into triangles and the spectral
problem for $L$ splits into separate spectral problems for each
triangle.

Let $\psi(x)=\psi(x;a_1,a_2,p_1,p_2)$ be a double-Bloch
eigenfunction for $L$, which corresponds to a point
$(a_1,a_2,p_1,p_2)$ of the surface \eqref{var} in accordance with
the formulas \eqref{ansatz}, \eqref{up}--\eqref{uppe} and
\eqref{p}. Given such a $\psi$, let us symmetrize it in the
following way:
\begin{equation}\label{spec}
\Psi(x)=\sum_{w\in W} \psi(wx)\,,
\end{equation}
where $W$ denotes the Weyl group for the system $B_2$. To get a
square-integrable eigenfunction, according to Proposition
\ref{sp1}, one takes $k\in 2\pi iP$, that is $k=(i\pi m, i\pi n)$
with integer $m,n$ having the same parity. For such $k$ the
substitution  of \eqref{p} into \eqref{var} leads to a system of
(transcendental) equations on $a_1,a_2$, and by solving it one
eventually finds the corresponding eigenfunctions. Note that
because of the invariance of the system \eqref{var} (and of
$\psi$) under the shifts \eqref{a}--\eqref{ak}, it is enough to
look for the solutions $(a_1,a_2)$ with $a_i$ lying inside the
fundamental parallelogram with the vertices $\pm (1+\tau)/2$.
Also, taking into account the Weyl group action, we can restrict
ourselves to the dominant weights, i.e. $0\le m\le n$.

Now let us consider the trigonometric limit $\tau\to +i\infty$,
then one can show that for any $m,n$ the corresponding system
\eqref{var} will have the unique solution $a_1,a_2$ inside the
fundamental parallelogram (at least, for sufficiently big $\tau$).
Moreover, if one fixes $k=(k_1,k_2)$ and takes then the limit
$\tau\to +\infty$, then $\psi$ will go to the Baker--Akhiezer
function $\psi(k,x)$ considered in \cite{C00} (this can be seen
directly from the formulas for $\psi$). Now, for the
Baker--Akhiezer function $\psi(k,x)$ it is known (see Theorem 6.7
of \cite{C00}) that the formula \eqref{spec} will produces all the
Jack polynomials if
\begin{equation}\label{adm}
k=2\pi i(\lambda+\rho)\qquad\text{with}\ \lambda\in P_+\quad\text{
and}\ \rho=\frac12\sum_{\alpha\in R_+} (m_\alpha+1)\alpha\,.
\end{equation}
For others $k\in 2\pi i P_+$ the symmetrized $\Psi$ will be zero.

In our situation this means that the function $\Psi=\Psi_{m,n}$
defined by \eqref{spec} will be non-zero as soon as
\begin{equation} \label{triangle} n-4\ge m\ge
2\,,\qquad \text{with}\quad n\equiv m \ (\mathrm {mod} 2)\,.
\end{equation}
For other $k\in 2\pi i P_+$ it will be zero in the trigonometric
limit. But the result of \cite{KT} cited in Section \ref{spe}
claims that the family of the eigenfunctions of $L$ is analytic in
$p=e^{\pi i\tau}$ and specializes to the Jack polynomials at
$p=0$. Hence, if $\Psi_{m,n}=0$ in the trigonometric limit, it
must be zero for all $\tau$ identically.

We conclude that the eigenfunctions $\Psi_{m,n}$ of $L$ are
labeled by $m,n$ satysfying \eqref{triangle}. In particular, the
ground state corresponds to $(m, n)=(2,6)$. The constructed
solutions $\Psi_{m,n}$ will have second order zeros along the
lines \eqref{sing} and will be invariant under orthogonal
reflections with respect to these lines. According to the results
of \cite{KT}, the resulting family is complete in $L^2(T)^W$ (see
Section \ref{spe} above).

\begin{remark}
For certain $(a_1,a_2,k_1,k_2)$ the corresponding Bloch solution
$\psi$ has a nontrivial symmetry in $x$, which must be a subgroup
of the Weyl group $W$. Our formulas for $\psi$ do not work
directly for some of these cases, because of the presence of an
extra component \eqref{vert}. Of particular interest are those of
the points which correspond to the solutions which are
double-(anti)periodic (in each of the variables $x_1,x_2$). These
can be viewed as multidimensional analogues of the classical
Lam\'e polynomials \cite{BE}.
\end{remark}

\section{Difference $B_2$ case} In this section we will
generalize the results above to the following difference version
of the operator \eqref{b2}:
\begin{equation}
\label{qb2} L= a_0 + a_+T_1^{2\omega}+ a_-T_1^{-2\omega}+
b_+T_2^{2\omega}+ b_-T_2^{-2\omega}\,,
\end{equation}
where $T_i^{\epsilon}$ stands for a shift in $x_i$ by $\epsilon$,
and the coefficients $a_\pm, b_\pm$ are
\begin{align*}
a_{\pm}=&\frac{\theta(x_1\mp\omega)\theta(x_1+x_2\mp 2\omega)
\theta(x_1-x_2\mp 2\omega)}{\theta(x_1\pm\o)\theta(x_1+x_2)
\theta(x_1-x_2)}\,,\\
b_{\pm}=&\frac{\theta(x_2\mp\omega)\theta(x_1+x_2\mp 2\omega)
\theta(x_1-x_2\pm 2\omega)}{\theta(x_2\pm\o)\theta(x_1+x_2)
\theta(x_1-x_2)}\,,
\end{align*}
while $a_0$ has the form $a_0=c_++c_-+d_++d_-$ with $c_\pm, d_\pm$
given by
\begin{align*}
c_{\pm}=&\frac{\t(2\o)\t(x_1\pm 5\o)\t(x_1+x_2\mp
2\o)\t(x_1-x_2\mp 2\o)}
{\t(4\o)\t(x_1\pm\o)\t(x_1+x_2)\t(x_1-x_2)}\,,\\
d_{\pm}=&\frac{\t(2\o)\t(x_2\pm 5\o)\t(x_1+x_2\mp
2\o)\t(x_1-x_2\pm 2\o)}
{\t(4\o)\t(x_2\pm\o)\t(x_1+x_2)\t(x_1-x_2)}\,.
\end{align*}
In all formulas $\theta(z)=\theta(z|\tau)$ is the odd Jacobi theta
function \eqref{theta1}.

This is a very special case of the so-called $BC_n$ generalization
of the quantum Ruijsenaars model \cite{R}. In trigonometric case
it has been introduced by Koornwinder \cite{Ko}. Elliptic version
was first suggested by van Diejen \cite{vD} and later extended in
\cite{Mi}, where its complete integrability has been proven. This
also can be viewed as an elliptic generalization of one of the
Macdonald operators \cite{Ma} for $B_2$. In what follows we assume
that the parameter $\omega$ is generic. Note that the operator
\eqref{b2} can be restored (up to a certain gauge) in the limit
$\omega\to 0$.

It is worth mentioning that coefficients of $L$ are not periodic,
so instead of (double-) Bloch solutions one should look for the
eigenfunctions in a certain $\t$-functional space. Another way of
putting it is to observe that $L$ can be reduced to elliptic form
using proper gauge. For instance, consider $$\widetilde
L=\delta^{-1}\circ L\circ\delta\,,\qquad
\delta=\t(x_1)\t(x_2)\t(x_1+x_2)\t(x_1-x_2)\,.$$ Then $\widetilde
L$ will have elliptic coefficients, so we can look for its Bloch
solutions $\psi(x)$. Correspondingly, $\Phi=\delta\psi$ will be an
eigenfunction for $L$ and it will have translation properties
similar to those of $\delta$. Abusing the language, below we refer
to $\Phi$ as a Bloch solution for $L$.

\subsection{Bloch solutions}
We are going to construct eigenfunctions of $L$ similar to the
differential case above. Our ansatz for $\Phi$ remains unchanged:
\begin{equation}\label{qphi}
\Phi=\exp(K_1x_1+K_2x_2)\sum_{0\le i,j\le2}
c_{ij}\theta\genfrac{[}{]}{0pt}{}{i/3}{0}(3x_1+\gamma_1|3\tau)
\theta\genfrac{[}{]}{0pt}{}{j/3}{0}(3x_2+\gamma_2|3\tau)\,.
\end{equation}
An analogue of the vanishing conditions \eqref{ax1}--\eqref{ax-}
is dictated by the singularities of $L$ and is the following:
\begin{align}\label{qax1}
\Phi(\omega,x_2) \equiv &\Phi(-\omega,x_2) && \quad \text{for
all}\ x_2\,,
\\ \label{qax2}
\Phi(x_1,\omega) \equiv &\Phi(x_1,-\omega) && \quad \text{for
all}\ x_1\,,
\\ \label{qax+}
  \Phi(x_1+\omega,x_2+\omega) \equiv &\Phi(x_1-\omega,x_2-\omega)
   && \quad \text{for}\ x_1+x_2=0\,, \\
  \label{qax-}
  \Phi(x_1+\omega,x_2-\omega) \equiv &\Phi(x_1-\omega,x_2+\omega)
   && \quad \text{for}\ x_1-x_2=0\,.
\end{align}
We are going to show that for a certain two-dimensional variety in
the space of parameters $\gamma_1,\gamma_2,K_1,K_2$ there is only
one (up to a factor) such $\Phi$. As a consequence, $\Phi$ will be
an eigenfunction of $L$, due to a natural analogue of Proposition
\ref{f}.

We start from conditions \eqref{qax+}, \eqref{qax-}. Using the
formula \eqref{add} and repeating the arguments used in case
$\omega=0$, we obtain a linear system for $c_{ij}$ and can see
that (for generic $\gamma_1,\gamma_2$) it defines $c_{ij}$
uniquely, up to a factor. However, solving this system leads to a
very cumbersome formula. Instead, let us define $\Phi$ by the
following formula:
\begin{equation}\label{qans}
\Phi=\exp(k_1x_1+k_2x_2) \sum_{i,j}b_{ij}\theta(x_1+a_1+i\omega)
\theta(x_2+a_2+j\omega)e^{\omega(ik_1+jk_2)}\,,
\end{equation}
where the summation is taken over the following set of indices:
$$(i,j)=(0,4),(4,0),(0,-4),(-4,0),(2,2),(2,-2),(-2,-2),(-2,2),(0,0)\,,$$
and the coefficients $b_{ij}=b_{ij}(x)$ look as follows:
$$b_{ij}=\beta_{ij}\theta\left(x_1+x_2-\frac{i+j}{2}\omega\right)
\theta\left(x_1-x_2- \frac{i-j}{2}\omega\right)$$ with
\begin{align*}
&\beta_{04}=\beta_{40}=\beta_{0,-4}=\beta_{-4,0}=(\theta(2\omega))^2
\\&\beta_{22}=\beta_{2,-2}=\beta_{-2,-2}=\beta_{-2,2}=-\theta(2\omega)
\theta(4\omega)\\ &\beta_{00}=(\theta(4\omega))^2\,.
\end{align*}

\begin{prop} For any $K_1,K_2$ and generic $\gamma_1,\gamma_2$ there
exists unique (up to a factor) function $\Phi$ of the form
\eqref{qphi} satisfying conditions \eqref{qax+}--\eqref{qax-}. It
is given by the formula \eqref{qans}, where $k_j=K_j+\pi i$ and
$a_j=\gamma_j-(1+\tau)/2$ ($j=1,2$).
\end{prop}

To prove the proposition one first checks that this $\Phi$ has the
needed translation properties in $x_1,x_2$, then an elementary
check shows that conditions \eqref{qax+}--\eqref{qax-} are
satisfied.

\medskip
Let us turn now to conditions \eqref{qax1}--\eqref{qax2}. First,
let us remark that \eqref{qax+}, \eqref{qax-} at $x=(0,0)$ imply
that
\begin{equation}\label{q0}
\Phi(\omega,\omega)-\Phi(-\omega,-\omega)=\Phi(\omega,-\omega)
-\Phi(-\omega,\omega)=0\,.
\end{equation}
Introduce now
\begin{gather*}
G_1=\Phi(\omega,\omega)-\Phi(\omega,-\omega)-\Phi(-\omega,\omega)
+\Phi(-\omega,-\omega)\,,\\
G_2=\Phi(\omega,-3\omega)-\Phi(-\omega,-3\omega)-\Phi(\omega,3\omega)
+\Phi(-\omega,3\omega)\,,
\end{gather*}
and consider the system
\begin{equation}\label{qob}
G_1=0\,,\quad G_2=0\,.
\end{equation}
Obviously, the vanishing condition \eqref{qax1} implies
\eqref{qob}. Conversely, from the system \eqref {qob} we deduce
immediately that the function
$$f(t)=\Phi(\omega,t)-\Phi(-\omega,t)$$ satisfies the conditions
$f(\omega)=f(-\omega)$ and $f(3\omega)=f(-3\omega)$. Together with
\eqref{q0} this gives that $f(\omega)=f(-\omega)=0$. Since $f$ is
a theta function of order $3$, it must have the form
$f=\theta(t-\omega)\theta(t+\omega)g(t)$ with
$g(3\omega)=g(-3\omega)$. Now, since $g$ is a theta function of
the first order with known characteristics (expressed in terms of
$a_2, k_2$), the condition $g(3\omega)=g(-3\omega)$ implies
$g\equiv 0$ as soon as
\begin{equation}\label{qgen}
e^{6\omega k_2}\ne
\frac{\theta(a_2-3\omega|\tau)}{\theta(a_2+3\omega|\tau)}\,.
\end{equation}
The latter condition is a difference version of \eqref{gen}. Now
let us collect some corollaries of \eqref{qax+}--\eqref{qax-}:
\begin{gather*}
\Phi(\omega,-3\omega)=\Phi(3\omega,-\omega)\,,\quad\Phi(-\omega,-3\omega)
=\Phi(-3\omega,-\omega)\,,\\\Phi(\omega,3\omega)=\Phi(3\omega,\omega)\,,
\quad \Phi(-\omega,3\omega)=\Phi(-3\omega,\omega)\,.
\end{gather*}
Thus, $G_2$ can be rewritten as $$G_2=\Phi(3\omega,-\omega)-
\Phi(-3\omega,-\omega)- \Phi(3\omega,\omega)
+\Phi(-3\omega,\omega)\,.$$ Hence, we can repeat the same
arguments with respect to $x_2$-variable and conclude that the
system \eqref{qob} implies \eqref{qax2} as soon as
\begin{equation}\label{qgen2}
e^{6\omega k_1}\ne
\frac{\theta(a_1-3\omega|\tau)}{\theta(a_1+3\omega|\tau)}\,.
\end{equation}
Summing up, we see that the system \eqref{qob} implies both of the
conditions \eqref{qax1}--\eqref{qax2} under assumptions
\eqref{qgen}--\eqref{qgen2}.

To get rid of restrictions \eqref{qgen}--\eqref{qgen2} let us
calculate $G_1, G_2$. First, introduce the notation $\xi_1,\xi_2$
for $$\xi_1=e^{\omega k_1}\,,\quad \xi_2=e^{\omega k_2}\,.$$ A
direct substitution gives that
\begin{align*}
\Phi(\o,\o)=&\beta_{0,-4}{\theta}(4\o){\theta}(-2\o)\t(a_1+\o)
\t(a_2-3\o)\xi_1\xi_2^{-3}+\\&
\beta_{-4,0}{\theta}(4\o){\theta}(2\o)\t(a_1-3\o)
\t(a_2+\o)\xi_1^{-3}\xi_2+\\&
\beta_{-2,2}{\theta}(2\o){\theta}(2\o)\t(a_1-\o)
\t(a_2+3\o)\xi_1^{-1}\xi_2^{3}+\\&
\beta_{2,-2}{\theta}(2\o){\theta}(-2\o)\t(a_1+3\o)
\t(a_2-\o)\xi_1^3\xi_2^{-1}\,.
\end{align*}
In a similar way we calculate $\Phi(-\o,\o)$ and
$G_1=2\Phi(\omega,\omega)-2\Phi(-\omega,\omega)$. As a result, the
equation $G_1=0$ takes the following form (up to a nonessential
factor):
\begin{multline}\label{aha1}
(\theta(a_1+3\omega)\xi_1^3-\theta(a_1-3\omega)\xi_1^{-3})
(\theta(a_2+\omega)\xi_2-\theta(a_2-\omega)\xi_2^{-1})\\=
(\theta(a_1+\omega)\xi_1-\theta(a_1-\omega)\xi_1^{-1})
(\theta(a_2+3\omega)\xi_2^3-\theta(a_2-3\omega)\xi_2^{-3})\,.
\end{multline}

In a similar way one can compute $G_2$. It turns out that it is a
linear combination of $8$ terms of the form $\xi_1^{\pm
3}\xi_2^{\pm 5}$, $\xi_1^{\pm 5}\xi_2^{\pm 3}$ and $8$ terms of
the form $\xi_1^{\pm 1}\xi_2^{\pm 3}$, $\xi_1^{\pm 3}\xi_2^{\pm
1}$. Moreover, the combination of the last $8$ terms is
proportional to $G_1$, so we can get rid of them by subtracting
$G_1$, this does not affect the system \eqref{qob}. After these
transformations equation $G_2=0$ takes the following nice form:
\begin{multline}\label{aha2}
(\theta(a_1+5\omega)\xi_1^5-\theta(a_1-5\omega)\xi_1^{-5})
(\theta(a_2+3\omega)\xi_2^3-\theta(a_2-3\omega)\xi_2^{-3})\\=
(\theta(a_1+3\omega)\xi_1^3-\theta(a_1-3\omega)\xi_1^{-3})
(\theta(a_2+5\omega)\xi_2^5-\theta(a_2-5\omega)\xi_2^{-5})\,.
\end{multline}

Summarizing, we see that the system \eqref{qob} is equivalent to
the equations \eqref{aha1}--\eqref{aha2}. These equations are
obviously invariant under the transformations
\eqref{a}--\eqref{ak}, in this way they define a covering over the
product $\mathcal E\times \mathcal E$ of two elliptic curves
$\mathcal E=\mathbb C/\mathbb Z + \tau \mathbb Z$. It has two
'vertical' components
\begin{equation}\label{qvert}
\{a_1=a_2\,,\ \xi_1=\xi_2\}\,,\qquad\{a_1=-a_2\,,\
\xi_1=(\xi_2)^{-1}\}\,.
\end{equation}
Another 'trivial' component is given by
\begin{equation}\label{qtriv}
\theta(a_1+3\omega)\xi_1^3-\theta(a_1-3\omega)\xi_1^{-3}=
\theta(a_2+3\omega)\xi_2^3-\theta(a_2-3\omega)\xi_2^{-3}=0\,.
\end{equation}
After deleting these three components, one gets a finite (in fact,
$17$-fold) covering of $\mathcal E\times \mathcal E$, let us
denote it by $\mathcal C$. We can conclude now that for any point
in $\mathcal C$ the corresponding function $\Phi$ will satisfy the
vanishing conditions \eqref{qax1}--\eqref{qax-} and, hence, it
will be an eigenfunction of the difference operator $L$.

\begin{theorem} Formulas \eqref{qans} and equations
\eqref{aha1}--\eqref{aha2} describe the (double-) Bloch
eigenfunctions of the difference operator \eqref{qb2}. The
Bloch--Hermite variety $\mathcal C$, which is obtained by deleting
components \eqref{qvert}--\eqref{qtriv} from the variety
\eqref{aha1}--\eqref{aha2}, is a $17$-fold covering of the product
$\mathcal E\times \mathcal E$ of two elliptic curves.
\end{theorem}

As a corollary, considering the limit $\omega\to 0$ we can
calculate explicitly the Hermite--Bloch variety for the operator
\eqref{b2}. Namely, one picks up the terms of order $4$ and $6$ in
$\omega$ in equation \eqref{aha1}.

\begin{cor} In the limit $\omega\to 0$ the system
\eqref{aha1}--\eqref{aha2} goes to \eqref{var}.
\end{cor}

The only thing we still have to explain is why the degree of the
covering is $17$. To this end, let us define a family of plane
rational curves $\varphi: \mathbb P^1\to\mathbb P^2$ of degree
$5$, depending on parameter $a\in \mathcal E$. Namely, for $a\in
\mathcal E$ and $u=(u_0:u_1)\in\mathbb P^1$ put
\begin{align*}
&\varphi(a,u)=(\varphi_0:\varphi_1:\varphi_2)\,,\quad\text{where}\\
&\varphi_0= \frac{ u_0^2u_1^3 \theta(a+\omega)\theta(a-{\o}/2)}
{\theta(a+{\o}/2)}- \frac{ u_0^3u_1^2
\theta(a-\omega)\theta(a+{\o}/2)} {\theta(a-{\o}/2)} \,,
\\&\varphi_1=
\frac{u_0u_1^4\theta(a+3\omega)\theta^3(a-{\o}/2)}
{\theta^3(a+{\o}/2)}-
\frac{u^4_0u_1\theta(a-3\omega)\theta^3(a+{\o}/2)}
{\theta^3(a-{\o}/2)}\,,
\\&\varphi_2=\frac{u^5_1\theta(a+5\omega)\theta^5(a-{\o}/2)}
{\theta^5(a+{\o}/2)}-
\frac{u^5_0\theta(a-5\omega)\theta^5(a+{\o}/2)}
{\theta^5(a-{\o}/2)}\,.
\end{align*}
Then the solutions $(\xi_1,\xi_2)$ of \eqref{var} correspond to
the intersection points of two curves $C_1=\varphi(a_1,\cdot\,)$,
$C_2=\varphi(a_2,\cdot\,)$ from our family. Namely, if
$\varphi(a_1,u)=\varphi(a_2,v)$ then $\xi_1, \xi_2$ with
$$(\xi_1)^2=\frac{u_1\theta^2(a_1-{\omega}/2)}
{u_0\theta^2(a_1+{\omega}/2)}\quad\text{and}\quad
(\xi_2)^2=\frac{v_1\theta^2(a_2-{\omega}/2)}
{v_0\theta^2(a_2+{\omega}/2)}$$ clearly satisfy
\eqref{aha1}--\eqref{aha2} and vice versa, provided
\eqref{qgen}--\eqref{qgen2}. We should, however, exclude from
consideration points with $\xi_1,\xi_2=0, \infty$ since
$\xi_i=e^{\o k_i}$. Namely, all the curves from our family pass
through $(0:0:1)=\varphi(a,0)=\varphi(a,\infty)$. A little
difference with $\o=0$ case is that now we have ${\rm
mult}\,(C_1\cap C_2)=8$ at this point. By Bezout's theorem, the
number of intersection points of $C_1, C_2$ , apart from
$(0:0:1)$, equals $5\times 5-8=17$.

\medskip
\begin{remark} Note that for given
$\xi_1^2, \xi_2^2$ the corresponding quasimomenta $k_1,k_2$ seem
to be non-unique, with the ambiguity of adding some multiple of
$i\pi/\o$. However, this would result in multiplying $\Phi$ by
$e^{i\pi x_1/\o}\,, e^{i\pi x_2/\o}$ which are quasi-constant on
the lattice $2\o\mathbb Z^2$. Thus, this leads to the same
eigenfunction, so the constructed Bloch solutions are in
one-to-one correspondence with the points of the surface $\mathcal
C$.
\end{remark}

\begin{remark}
It is clear that the Weyl group action on $\mathcal C$ is
generated by two involutions $$(a_1,\xi_1)\to
(-a_1,\xi_1^{-1})\quad\text{and }(a_1,\xi_1)\leftrightarrow
(a_2,\xi_2)\,.$$
\end{remark}

\subsection{Structure of the solution space}

As we mentioned above, the difference operator \eqref{qb2} is
completely integrable, i.e. there exists another difference
operator $L_1$ which commutes with $L$. It is given by the
following expression \cite{vD,Mi}:
\begin{equation}\label{qL1}
L_1=c_{++}T^{\o}_1T^{\o}_2+c_{+-}T^{\o}_1T^{-\o}_2+c_{-+}T^{-\o}_1T^{\o}_2+
c_{--}T^{-\o}_1T^{-\o}_2\,,
\end{equation}
where the coefficients $c_{\epsilon_1,\epsilon_2}$ look as follows
(we treat $\pm$ as $\pm 1$):
\begin{equation*}
c_{\epsilon_1,\epsilon_2}=\frac
{\t(x_1-\o\epsilon_1)\t(x_2-\o\epsilon_2)
\t(x_1+x_2-2\o(\epsilon_1+\epsilon_2))\t(x_1-x_2-2\o(\epsilon_1-\epsilon_2))}
{\t(x_1)\t(x_2)\t(x_1+x_2)\t(x_1-x_2)}\,.
\end{equation*}

\medskip
Now let us consider the system of two partial difference
equations:
\begin{equation}\label{qs}
Lf=Ef\,,\quad L_1f=E_1f\,,
\end{equation}
defined on the lattice $2\o\mathcal L$ where $$\mathcal
L=\{(m,n)\mid m\pm n\in\mathbb Z\}\,.$$ More precisely, we fix
generic $x^0\in\mathbb C^2$ as a base point and regard a function
$f$ in \eqref{qs} as being defined on $x^0+2\o\mathcal L\subset
\mathbb C^2$. The base point $x^0$ must be outside the singular
locus of $L, L_1$, i.e. such that $L, L_1$ are nonsingular on
$x^0+2\o\mathcal L$.

The Bloch solutions $\Phi$, constructed above, are common
eigenfunctions of $L$ and $L_1$. (The proof for $L_1$ is the same:
the only thing to check is an analogue of Proposition \ref{f}.)
Since $L$ and $L_1$ are $W$-symmetric, each of the $8$ functions
$\Phi(wx)$ ($w\in W$) will solve the system \eqref{qs}. We know
that for generic point of the spectral surface $X$ (thus, for
generic $E, E_1$) all $8$ functions $\Phi(wx)$ are linearly
independent (as functions on  $\mathbb C^2$), because they have
different translation properties with respect to the shifts
\eqref{a}--\eqref{ak}. Hence, their restriction to
$x^0+2\o\mathcal L$ also gives $8$ linearly independent solutions
of \eqref{qs}(at least, for generic base point $x^0$). On the
other hand, it is not difficult to see that any solution $f$ is
uniquely determined by its values at eight points $x^0+\nu$ with
the following $\nu$: $$\nu=(0,0)\,,(\pm\o,\o)\,,(\o,-\o)\,,(\pm
2\o,0)\,,(0,2\o)\,,(\o,3\o)\,.$$ This implies that the dimension
of the solution space of \eqref{qs} is at most $8$. Thus, we
conclude that any solution of \eqref{qs} (for generic $E, E_1$) is
a linear combination of $8$ Bloch solutions $\{\Phi(wx)\}_{w\in
W}$.

\begin{prop} The space of solutions of the system \eqref{qs} has
dimension $8$ and for generic $E, E_1$ is generated by the Bloch
solutions $\Phi(wx)$ ($w\in W$).
\end{prop}

\begin{remark} Above we associated a double-Bloch eigenfunction $\Phi$ to
a solution $(a_1,a_2,\xi_1,\xi_2)$ of the equations \eqref{aha1}--
\eqref{aha2}. Note that these equations are invariant under
$\xi_j\to -\xi_j$, but this does not lead to another
eigenfunction, since they will differ by a quasiconstant factor.
Situation is different for the system \eqref{qs}, since it is
defined on a different lattice. It is easy to see that
$(\xi_1,\xi_2)$ and $(-\xi_1,-\xi_2)$ still lead to the same
solution modulo quasiconstants, the same is true for
$(\xi_1,-\xi_2)$ and $(-\xi_1,\xi_2)$. The resulting two functions
have the same eigenvalue $E$ in \eqref{qs}, but opposite values of
$E_1$. Thus, the Bloch variety for the system \eqref{qs} is a
double covering of the surface $\mathcal C$ introduced above.
\end{remark}

\medskip
We will not go into discussing the spectral properties of the
difference operator $L$. See papers \cite{R1,Kom} devoted to this
rather delicate matter. Let us just remark on some special
solutions analogous to the 'discrete spectrum' considered in
section \ref{spectrum}. Namely, let us consider the following
anti-invariant solution of the system \eqref{qs}:
\begin{equation*}
\Phi_{skew}(x)=\sum_{w\in W}({\det w})\Phi(wx)\,.
\end{equation*}
The vanishing conditions \eqref{qax1}--\eqref{qax-} imply that
$\Phi_{skew}$ vanishes along lines $x_1\pm x_2=\pm2\o$ and
$x_j=\pm\o$. It also vanishes if $x_1=\pm x_2$ and $x_j=0$ due to
anti-invariance. Let us require now for all $8$ functions
$\Phi(wx)$ to have the same Floquet--Bloch multipliers with
respect to the shifts by $e_1$ and $e_2$, which is equivalent to
the conditions
\begin{equation*}
\exp (k_1)=\exp (k_2)=\pm 1\,.
\end{equation*}
Then $\Phi_{skew}$ will vanish also along the shifted lines
$$x_1\pm x_2=m\,, m\pm2\o\,,\qquad x_j=n\,, n\pm\o\quad
(m,n\in\mathbb Z)\,.$$ In the limit $\o\to 0$ these solutions  go
to those $\Psi$ constructed in section \ref{spectrum}, more
precisely, $$\o^{-6}\Phi_{skew} \longrightarrow 64
(\t'(0))^2\t(a_1)\t(a_2)
{\theta}(x_1){\theta}(x_2){\theta}(x_1+x_2){\theta}(x_1-x_2)
\Psi\,\quad \text{as}\quad \o\to 0\,.$$

\section{Hietarinta operator and its discretization}

\subsection{Continuous case}
We consider now the Schr\"odinger operator \eqref{3} but first let
us rescale the coordinates $x_i\to a_ix_i$, so instead of
\eqref{3} we will consider
\begin{multline}\label{blL}
L=-a_1^2\d_1^2-a_2^2\d_2^2-a_3^2\d_3^2\\+2(a_1^2+a_2^2)\wp(x_1-x_2)+
 2(a_2^2+a_3^3)\wp(x_2-x_3)+2(a_3^2+a_1^2)\wp(x_3-x_1)\,,
\end{multline}
where, as before, $\wp(z)=\wp(z|\tau)$ is the Weierstrass
$\wp$-function and $a_1^2+a_2^2+a_3^2=0$. We are going to
calculate the double-Bloch eigenfunctions of $L$. More
specifically, we are looking for the solutions $\psi$ of the
equation $L\psi=E\psi$ with the following properties:

(i) $\psi$ is of the form
\begin{equation}\label{blan}
\psi(x)=\frac{\Phi(x)}{\theta(x_{12})\theta(x_{23})\theta(x_{31})}
\exp(k_1x_1+k_2x_2+k_3x_3)\,,
\end{equation}
where $x_{ij}:=x_i-x_j$,\
$\theta=\theta\genfrac{[}{]}{0pt}{}{1/2}{1/2}$ and $\Phi$ is
holomorphic in $\mathbb C^3$ and depends on the differences
$x_{ij}$ only, in other words, $(\d_1+\d_2+\d_3)\Phi=0$\,;

(ii) $\psi$ has the following translation properties:
\begin{equation}\label{re}
\psi(x+e_j)=e^{k_j}\psi(x)\,,\qquad\psi(x+\tau
e_j)=e^{\mu_j}\psi(x)\quad (j=1,2,3)\,,
\end{equation}
where $(e_1, e_2, e_3)$ is the standard basis in $\mathbb C^3$.

It is not difficult to conclude that for fixed $k_j,\mu_j$ the
conditions above determine a three-dimensional functional space,
and the corresponding $\Phi(x)$ in \eqref{blan} must be of the
form
\begin{equation}\label{blphi}
\Phi=\sum_{l=0}^2c_l\theta(x_{12}+b_{12}+l\tau/3)
\theta(x_{23}+b_{23}+l\tau/3) \theta(x_{31}+b_{31}+l\tau/3)\,,
\end{equation}
where $c_0,c_1,c_2$ are arbitrary constants and the parameters
$b_{12}, b_{23}, b_{31}$ are related to $\mu_j$ above in the
following way:
\begin{gather}\notag
e^{\mu_1}=e^{k_1\tau+2\pi ib_{31}-2\pi ib_{12}}\,,\\\label{pa}
e^{\mu_2}=e^{k_2\tau+2\pi ib_{12}-2\pi ib_{23}}\,,\\\notag
e^{\mu_3}=e^{k_3\tau+2\pi ib_{23}-2\pi ib_{31}}\,.
\end{gather}
This shows that the three-dimensional space \eqref{blphi} depends,
essentially, on the pairwise differences of the parameters
$b_{lm}$ only. Thus, without loss of generality we may assume that
\begin{equation}\label{zero}
b_{12}+b_{23}+b_{31}=0\,.
\end{equation}
In formulas below we will also use $b_{21},b_{32},b_{13}$ under
the convention that $ b_{ij}=-b_{ji}$.

Now, in accordance with Proposition \ref{f}, we impose certain
vanishing conditions on $\psi$ which are motivated by the
structure of the singularities of the operator \eqref{blL}.
Namely, for any $i=1,2,3$ consider the function
$f(t)=\psi(x+ta_{i-1}^2e_{i-1}-ta_{i+1}^2e_{i+1})$ for $x$ such
that $x_{i-1}=x_{i+1}$ (we treat indices modulo $3$, so
$x_0=x_3$). Our assumptions about $\psi$ imply that for such $x$
the function $f$ will have a pole at $t=0$, so its Laurent
expansion will look as $f=a_{-1}t^{-1}+a_0+a_1t+\dots$. The
coefficients in this expansion depend on $x$. Let us require that
$a_0=0$ for all $x$ such that $x_{i-1}=x_{i+1}$. Using
\eqref{blan} one rewrites this condition as follows:
 \begin{align}\label{nax1}
&a_{i-1}^2k_{i-1}-a_{i+1}^2k_{i+1}+ F_{i}=0\,,\\\label{Fi}
&F_{i}:=\frac{a_{i-1}^2
 \d_{i-1}\Phi-a_{i+1}^2\d_{i+1}\Phi}
 {\Phi}-a_{i-1}^2\frac{\theta'(x_{i-1,i})}{\theta(x_{i-1,i})}+
 a_{i+1}^2\frac{\theta'(x_{i+1,i})}{\theta(x_{i+1,i})}\,,
\end{align}
with \eqref{nax1} to be valid for all $x$ such that
$x_{i-1}=x_{i+1}$.

The following lemma follows from Proposition \ref{f}.
\begin{lemma}\label{blinv}
If $\psi$ has the form \eqref{blan}, \eqref{blphi} and satisfies
the vanishing conditions \eqref{nax1}, then the same will be true
for its image $\widetilde\psi=L\psi$ under the action of the
operator \eqref{blL}.
\end{lemma}

As we will see below, for a certain three-dimensional subvariety
in the space of the parameters $k_j,b_{lm}$ the vanishing
conditions cut a one-dimensional subspace in the space
\eqref{blphi}. Thus, the lemma ensures that the corresponding
$\psi(x)$ will be an eigenfunction of $L$.

We may regard the restriction of the expression $F_{i}$ on the
plane $x_{i-1}=x_{i+1}$ as a function of $z=x_{i-1,i}=x_{i+1,i}$.
It is easy to check then that $F_i$ will be an elliptic function
of $z$ with periods $1, \tau$. So, first of all we have to choose
the parameters $b_{lm}, c_j$ in such a way that $F_i(z)$ would be
non-singular. Let us assume that the parameters $a_1,a_2,a_3$ are
generic enough, i.e. that $a_i^2\ne a_j^2$. Then for $F_i$ to be
non-singular at $z=0$ we need $\Phi(0)=0$. This gives the
following condition:
 \begin{equation}\label{fcon}
 \widetilde c_0+\widetilde c_1+\widetilde c_2=0\,,\qquad
 \widetilde c_l=c_l\theta(b_{12}+\frac{l\tau}3)
 \theta(b_{23}+\frac{l\tau}3)
 \theta(b_{31}+\frac{l\tau}3)\,.
 \end{equation}
Now we note that $$\Phi(z):=\Phi\mid_{x_{i-1}=x_{i+1}}$$ is a one-
dimensional $\theta$-function of order $2$, hence it has two zeros
(modulo $1,\tau$). First zero is $z=0$ (due to condition
\eqref{fcon}). An easy check shows that the second zero is
$z=b_{i,i-1}+b_{i,i+1}$. So, up to a constant factor,
$$\Phi(z)=\theta(z)\theta(z-b_{i,i-1}-b_{i,i+1})\,.$$ To get rid
of a possible pole at $z=b_{i,i-1}+b_{i,i+1}$ in \eqref{Fi} we
must require that $a_{i-1}^2
 \d_{i-1}\Phi-a_{i+1}^2\d_{i+1}\Phi=0$ for
 $x_{i-1,i}=x_{i+1,i}=b_{i,i-1}+b_{i,i+1}$.
This leads to the following relation:
\begin{align}\label{scon}
\sum_{l=0}^2\widetilde c_l(a_1^2\zeta(b_{23}+\frac{l\tau}3)
+a_2^2\zeta(b_{31}+\frac{l\tau}3)
+a_3^2\zeta(b_{12}+\frac{l\tau}3))=0\,.
\end{align}
Here and below $\zeta(z):=\frac{d}{dz}log\theta(z)$. Notice that
the relation \eqref{scon} is symmetric with respect to indices
$1,2,3$ (so we have just one condition instead of possible
three!).

We use the relations \eqref{fcon}, \eqref{scon} to express (up to
a common factor) the parameters $c_l$ in terms of the parameters
$b_{ij}$. These relations imply that each of $F_i$ \eqref{Fi} is
nonsingular in $z=x_{i-1,i}=x_{i+1,i}$, therefore they are some
constants depending on $c_l, b_{ij}$. Thus, \eqref{nax1} leads to
the expressions for the differences
$a_{i-1}^2k_{i-1}-a_{i+1}^2k_{i+1}$ in terms of $c_l$ and
$b_{lm}$. One can check that the resulting system is always
compatible (i.e. that $F_1+F_2+F_3=0$). This follows, for
instance, from  the compatibility of the system  \eqref{qscon}
below by going to the limit $\o\to 0$. As a corollary, the
formulas \eqref{blan},\eqref{blphi}, together with \eqref{fcon},
\eqref{scon} and \eqref{nax1} deliver the expression for the
double-Bloch eigenfunctions of the operator $L$.

Finally, let us discuss the structure of the Hermite--Bloch
variety of the operator \eqref{blL}. The double-Bloch solutions
are parametrized by $b_{12}, b_{23}, b_{31}$ with $b_{12}+ b_{23}+
b_{31}=0$, and the corresponding $k_1,k_2,k_3$ are determined from
\eqref{nax1}. Denote by $b$ and $k$ the three-component vectors
$b=(b_{12}, b_{23}, b_{31})$ and $k=(k_1, k_2, k_3)$. Then two
different points in
 the parameter space $(b,k)$ lead
to the same solution iff the corresponding Floquet multipliers in
\eqref{re} are the same. Taking into account relations \eqref{pa},
we conclude that the following transformations do not lead to
another Bloch solution:
\begin{align}\notag
&b \to b+\varepsilon_1\,, \quad k\to k
\qquad(\varepsilon_1=(2/3,-1/3,-1/3))\\\label{bt} &b\to
b+\varepsilon_2\,,\quad k\to k
\qquad(\varepsilon_2=(-1/3,2/3,-1/3))
\\\notag &b\to b+\tau\varepsilon_1\,,
\quad k\to k+2\pi i(1,-1,0)\\\notag &b\to
b+\tau\varepsilon_2\,,\quad k\to k+2\pi i(0,1,-1)\,.
\end{align}
Thus, $b$ is effectively represented by a point of a factor
$\mathbb C^2/\mathcal L+\tau \mathcal L$ where $\mathbb
C^2=\{z_1+z_2+z_3=0\}\subset \mathbb C^3$ and the lattice
$\mathcal L$ is generated by $\varepsilon_1, \varepsilon_2$. This
factor is isomorphic to the product of two elliptic curves with
parameter $\tau$. Above each point $b$ we have a complex line of
double-Bloch solutions, because equations \eqref{nax1} determine
$k$ up to adding any multiple of $(a_1^{-2},a_2^{-2},a_3^{-2})$.

 \subsection{Discrete case}
 We keep the notation  $x_{ij}$ for $x_i-x_j$. The discrete
 version of the operator \eqref{blL} looks as follows:
\begin{equation}\label{qblL}
D=\sum_{i=1}^3\frac { \theta(\omega)\theta(x_{i-1,i}+\omega
a_i^2)\theta(x_{i,i+1}-\omega a_i^2)} {
 \theta(\omega a_i^2)
 \theta(x_{i-1,i})\theta(x_{i,i+1})
} T_i^{\omega a_i^2}\,,
\end{equation}
where $a_1^2+a_2^2+a_3^2=0$ and $T_i^\epsilon$ stands for a shift
by $\epsilon$ in $x_i$. Its rational version ${\theta}(z)=z$ was
communicated to us by M.Feigin who found it to be dual (in
bispectral sense) to the trigonometric version $\wp=\sin^{-2}$ of
the Hietarinta operator.

The difference operator \eqref{qblL} relates to \eqref{blL} in the
following way:
$$D=a_1^{-2}+a_2^{-2}+a_3^{-2}+\omega(\d_1+\d_2+\d_3)+\frac{\omega^2}{2}
({\rm const}-\widetilde L)+o(\omega^2) \quad\text{as}\quad
\omega\to 0\,,$$ where $\widetilde L$ is gauge-equivalent to $L$,
$$\widetilde L=\delta\circ L \circ \delta^{-1}\,,\qquad
\delta=\theta(x_1-x_2)\theta(x_2-x_3)\theta(x_3-x_1)\,.$$

Unlike $L$, the operator $D$ is not periodic. As a result, instead
of the  double Bloch eigenfunctions, we will look for
eigenfunctions with the translation properties similar to those of
$\delta$. Apart from that, our ansatz for the eigenfunctions
$\varphi$ of the operator $D$ remains the same:
\begin{align}\label{qblan}
&\varphi=\exp(k_1x_1+k_2x_2+k_3x_3)\Phi\,,\\\label{qblphi}
&\Phi=\sum_{l=0}^2c_l\theta(x_{12}+b_{12}+l\tau/3)
 \theta(x_{23}+b_{23}+l\tau/3)
 \theta(x_{31}+b_{31}+l\tau/3)\,,\\\label{qbb}
 &b_{12}+b_{23}+b_{31}=0\,.
\end{align}
The vanishing conditions now look as follows: for each $i=1,2,3$
\begin{equation}\label{BLax}
 F_i:=\theta(x_{i,i-1}+\omega a_{i-1}^2)
 T_{i-1}^{\omega a_{i-1}^2}(\varphi)-
 \theta(x_{i,i+1}+\omega a_{i+1}^2)
 T_{i+1}^{\omega a_{i+1}^2}(\varphi)=0
 \end{equation}
identically for all $x$ with $x_{i+1}=x_{i-1}$.

We have then a straightforward analog of Lemma \ref{blinv}, so the
same approach as above will give us the eigenfunctions for $D$.

Let us first formulate the result. Namely, we consider the
following two conditions on the function $\Phi(x_1,x_2,x_3)$ given
by \eqref{qblphi}:
\begin{equation}\label{qfcon}
\Phi(\omega a_1^2,0,-\omega a_2^2)=0\,,\qquad \Phi(-\omega
a_2^2,0,\omega a_3^2)=0\,.
\end{equation}
This gives us two linear equations on $c_l$ and we use them to
express $c_l$ (up to a factor) through $b_{ij}$.

Secondly, we impose the following three relations:
\begin{align}\notag
e^{\omega a_{1}^2k_{1}-\omega a_{3}^2k_{3}}=&\frac{\theta(\omega
a_{3}^2)}{\theta(\omega a_{1}^2)}\frac{\Phi(0,0,\omega
a_{3}^2)}{\Phi(\omega a_{1}^2,0,0)}\,\\\label{qscon} e^{\omega
a_{2}^2k_{2}-\omega a_{1}^2k_{1}}=&\frac{\theta(\omega
a_{1}^2)}{\theta(\omega a_{2}^2)}\frac{\Phi(\omega
a_{1}^2,0,0)}{\Phi(0,\omega a_{2}^2,0)}\,\\\notag e^{\omega
a_{3}^2k_{3}-\omega a_{2}^2k_{2}}=&\frac{\theta(\omega
a_{2}^2)}{\theta(\omega a_{3}^2)}\frac{\Phi(0,\omega
a_{2}^2,0)}{\Phi(0,0,\omega a_{3}^2)}\,.
\end{align}
We use these formulas to express $k_1,k_2,k_3$ through $\Phi$. The
solution is not unique, and in fact we have a one-parameter family
of $k_i$. Altogether, formulas \eqref{qfcon}--\eqref{qscon} fix
the dependence of $c_l$ and $k_i$ and hence of $\varphi$ on three
parameters $b_{ij}$ (related by \eqref{qbb}). The resulting family
of functions $\varphi(x)$ depends on three parameters: two of
$b_{ij}$ and one more due to the freedom in resolving
\eqref{qscon}, see more comments below.

\begin{theorem}
The formulas \eqref{qblan}--\eqref{qbb} and the relations
\eqref{qfcon}--\eqref{qscon} give a three-parameter family of
eigenfunctions for the difference operator $D$.
\end{theorem}

To prove the theorem, let us first notice that each of $F_i$ in
\eqref{BLax}, being regarded as a function of
$z=x_{i,i-1}=x_{i,i+1}$, is a one-dimensional theta-function of
order $3$, so if it doesn't vanish, it must have three zeros
(modulo $1, \tau$). Moreover, a simple count shows that the sum of
these zeros will be equal to $b_{i,i-1}+b_{i,i+1}$. On the other
hand, a direct substitution into \eqref{BLax} shows that the
relations \eqref{qfcon} imply that $F_2$ vanishes for $z=-\omega
a_1^2$ and $z=-\omega a_3^2$. Further, the first relation in
\eqref{qscon} simply encodes the fact that $F_2$ vanishes at
$z=0$. Since the sum of these three zeros is $\omega a_2^2$ which,
generically, is not $b_{21}+b_{23}$, we conclude that $F_2(z)$ is
zero identically.

At first glance it seems that we need to add four more conditions
to ensure all the vanishing properties \eqref{BLax}. Namely, one
needs also
\begin{multline}\label{more}
\Phi(0,\omega a_2^2,-\omega a_1^2)=\Phi(0,-\omega a_1^2,\omega
a_3^2)\\=\Phi(-\omega a_3^2,\omega a_2^2,0)=\Phi(\omega
a_1^2,-\omega a_3^2,0)=0\,.
\end{multline}
However, since $\Phi$ depends on the pairwise differences of $x_i$
only, we will have that $$\Phi(0,\omega a_2^2,-\omega
a_1^2)=\Phi(-\omega a_2^2,-\omega a_2^2+\omega a_2^2,-\omega
a_2^2-\omega a_1^2)=\Phi(-\omega a_2^2,0,\omega a_3^2)=0\,.$$ In
the same way other relations in \eqref{more} follow from
\eqref{qfcon}.

This demonstrates that the three-parameter family constructed in
the theorem satisfies the vanishing conditions \eqref{BLax}, thus
proving the theorem.

Finally, let us comment on the structure of the Hermite--Bloch
variety. Similarly to the case $\o=0$ above, Bloch solutions are
parametrized by $b=(b_{12}, b_{23}, b_{31})$ with $b_{12}+ b_{23}+
b_{31}=0$. This determines the corresponding $c_l$ by
\eqref{qfcon}. After that $k=(k_1,k_2,k_3)$ are determined from
\eqref{qscon}. At this point we have certain freedom: if
$k=(k_1,k_2,k_3)$ is a solution of \eqref{qscon}, then any
$$k'=k+\frac{t}{\o}(a_1^{-2}, a_2^{-2}, a_3^{-2})+\frac{2\pi
i}{\o} (n_1a_1^{-2}, n_2a_2^{-2}, n_3a_3^{-2})$$ with any
$t\in\mathbb C$ and integer $n_1,n_2,n_3$ will be a solution, too.
However, the last term is not essential since it results in
multiplying $\Phi$ by a quasiconstant. For the same reason, the
factor $t$ in the second term is essential modulo $2\pi i$ only.
Besides, we still have the translation invariance of $\Phi$ with
respect to the transformations \eqref{bt}. Thus, the
Hermite--Bloch variety is fibered over the product of two elliptic
curves with the fibers isomorphic to $\mathbb C/2\pi i$.

\end{document}